\documentclass[draft,preprint,12pt,numbers,sort&compress]{elsarticle}
\usepackage{}
\usepackage{mathrsfs}
\usepackage{amssymb}
\usepackage{amsthm}
\usepackage{mathrsfs}
\usepackage[centertags]{amsmath}
\usepackage{amsfonts}

\usepackage{color}

\input amssym.def
\input amssym.tex

\date{}

\newtheorem{Theorem}{Theorem}[section]

\newtheorem{Lemma}{Lemma}[section]

\newcommand\R{\mbox{\bf R}}

\newcommand\SR{\mbox{\scriptsize\bf R}}

\newcommand{\definition}{{\lower .5ex
  \hbox{$\>\>\stackrel{\triangle}{=}\>\>$} }}
\newcommand\supp{\mathop{\rm supp}}

\begin{document}

\baselineskip=22pt
\thispagestyle{empty}

\mbox{}
\bigskip

\begin{center}
{\Large \bf Random  data Cauchy problem for a generalized    }\\[1ex]
{\Large \bf KdV equation  in the supercritical case}

{Wei Yan\footnote{ Email:yanwei19821115@sina.cn}$^{a,b}$, \quad  Jinqiao
Duan\footnote{Email:duan@iit.edu}$^b$,\quad Jianhua  Huang\footnote{Email:jhhuang32@nudt.edu.cn}$^c$}\\[1ex]
{$^a$College of Mathematics and Information Science, Henan Normal University,}\\
{Xinxiang, Henan 453007,   China}\\[1ex]

{$^b$ Department of Applied Mathematics, Illinois Institute of Technology,}\\[1ex]
{Chicago, IL 60616, USA  }\\[1ex]
{$^c$College of Science, National University of Defense  Technology,}\\
{ Changsha, Hunan 410073,  China}\\[1ex]

\end{center}

\bigskip
\bigskip

\noindent{\bf Abstract.}   We consider the Cauchy problem for a generalized
KdV  equation
\begin{eqnarray*}
      u_{t}+\partial_{x}^{3}u+u^{7}u_{x}=0,
    \end{eqnarray*}
 with random data on $\R$. Kenig, Ponce, Vega(Comm. Pure Appl. Math.46(1993), 527-620)
 proved that the  problem
is globally well-posed in $H^{s}(\R)$ with $s> s_{crit}=\frac{3}{14}$, which is the scaling
critical regularity indices.
Birnir,   Kenig,  Ponce, Svanstedt,  Vega(J. London Math. Soc. 53 (1996),
551-559.) proved that the problem is ill-posed
in the sense that the time of existence $T$ and the continuous dependence cannot
be expressed in terms of the size of the data in the $H^{\frac{3}{14}}$-norm.
In this present paper,  we prove that almost sure local in time
well-posedness
holds in $H^{s}(\R)$ with $s>\frac{17}{112},$ whose lower bound is below  $\frac{3}{14}.$
 The key ingredients are  the Wiener randomization of the initial data and
probabilistic  Strichartz estimates together with some important embedding Theorems.

\noindent {\bf Keywords}: Generalized KdV  equation; Almost sure local well-posedness; Modulation space;
Wiener decomposition;

\leftskip 0 true cm \rightskip 0 true cm

\newpage{}

\bigskip

{\large\bf 1. Introduction}
\bigskip

\setcounter{Theorem}{0} \setcounter{Lemma}{0}

\setcounter{section}{1}

In this paper, we consider the Cauchy problem for the generalized Korteweg-de
Vries(gKdV) equation
\begin{eqnarray}
      u_{t}+\partial_{x}^{3}u+u^{7}u_{x}=0, \label{1.01}
    \end{eqnarray}
with random data on $\R$.

By using Strichartz estimates and commutator estimates as well as interpolation theorems,
Kenig, Ponce, Vega \cite{KPV1993} proved that the  Cauchy problem for (\ref{1.01})
is globally well-posed in $H^{s}(\R)$ with $s> s_{crit}=\frac{3}{14}$.
Obviously, if $u(x,t)$  is     the  solution to   (\ref{1.01}),  then $u_{\lambda}(x,t)=\lambda^{-\frac{2}{7}}u\left(\frac{x}{\lambda},\frac{t}{\lambda^{3}}\right)$ is the solution to
 (\ref{1.01})  with the $\lambda$-scaled
initial data $u_{0\lambda}=\lambda^{-\frac{2}{7}}u_{0}\left(\frac{x}{\lambda}\right)$.
When $s_{crit}=\frac{3}{14}$, it is easily checked that $\|u_{0\lambda}(x)\|_{\dot{H}^{s_{crit}}}=
\|u_{0}(x)\|_{\dot{H}^{s_{crit}}}$. When $s<s_{crit}$, this is the supercritical regime. Birnir  et al.
 \cite{BKPSV} proved that the problem is ill-posed in $H^{\frac{3}{14}}(\R)$
 in the sense that the time of existence $T$ and the continuous dependence cannot
be expressed in terms of the size of the data in the $H^{\frac{3}{14}}$-norm. By using homogeneous and
nonhomogeneous Besov spaces,
Molinet and  Ribaud \cite{MR}   considered the local and global Cauchy
problem for the k-gKdV equation and  existence and uniqueness of similarity solutions.  Merle
\cite{M2001} studied the  existence of
 blow-up solutions in the energy space for the
critical generalized gKdV equation.
 Martel and Merle \cite{MM} studied the  instability of solitons for the critical gKdV equation.
Martel and  Merle \cite{MMJAMS}  studied the  blow up in
finite time and dynamics of blow
up solutions for the $L^{2}$-critical gKdV equation.
Martel and  Merle \cite{MM2002}  studied the stability of blow-up profile and lower bounds
for blow-up rate for the critical gKdV equation. By using the
smoothing effect of Kato's type associated to the linear problem, the
 fractional derivative commutators and chain rule for fractional derivatives
  as well as some interpolation theorems,
Fonseca et al.
 \cite{FLP}  proved the global well-posedness of  quintic gKdV  equation in
 $H^{s}(\R)$ in $s > 3/5$
assuming $\|u_{0}\|_{L^{2}}\leq \|Q\|_{L^{2}}$ with $Q=[3c sech^{2}(2\sqrt{c}x)]^{1/4}.$
By using the $X_{s,b}$ spaces,
Colliander et al. 	\cite{CKSTT}  studied the periodic  quartic and quintic
 gKdV equation on the torus. By using the $I$-method,
 Farah \cite{F}  proved that
the initial value problem for the critical gKdV on the real line is globally
well-posed in $H^{s}(\R)$ in $s > 3/5$ with the appropriate smallness assumption
on the initial data. By using the Strichartz estimates and Airy linear profile decomposition,
  Killip et al. \cite{KKSV} studied the blow-up of  the mass-critical gKdV.
   By using the I-method and
   the multilinear correction analysis,
Miao  et al. \cite{MSWX} proved that
the initial value problem for the critical gKdV  is globally
well-posed in $H^{s}(\R)$ in $s > \frac{3}{16}.$
Recently, by using Strichartz estimates and commutator estimates,  Kenig et al. \cite{KLPV}
studied the special regularity properties of solutions to
 the k-gKdV equation.

 Building upon \cite{LRS} and  using the $X_{s,b}$-spaces defined below and Strichartz estimates
and probabilistic tools,  Bourgain \cite{B1994,B-1994, B-1996} extended the local solution to
global solution with large set of initial data. The approach initiated by Bourgain
\cite{B1994,B-1994, B-1996}
attracted the attention of many people around the
problem of constructing the invariant measures for
many evolution equations with random initial data and  constructing large set of initial data  of
supercritical regime
 \cite{BT2007, B-2014, BB2014, BOP,BO-tms,BT2008-I,BT2008,BTT,BT2014,CG,CO,D,DC,FS,HO,HO15,
 LM,LRS,NO,NPS,O2009,Oh,OP,OQV,OQ2013,ORT,OSIAM,P,Poiret12, Poiret1207,PRT,R,ST,TT,T,T2009, TV2013, TV, ZFPAMS, ZF}
  and the references therein.
After a suitable randomization, by using probabilistic Strichartz estimates,  Burq and Tzvetkov
\cite{BT2008-I,BT2008} constructed
local and global  strong solution for a large set of initial data for the supercritical wave
equation on  three dimensional
compact Riemannian manifold. It is by now well understood that probabilistic tools play an
important role in extending the local solution to global solution
with  large set of initial data and constructing the solution in the supercritical regime
with a large set of initial data.
In the absence of  invariant measures, by suitably adapting Bourgain's high-low frequency
decomposition,   Colliander and Oh \cite{CO} studied almost sure global well-posedness for
the subcritical 1D periodic cubic  nonlinear
 cubic Schr\"odinger equation  below $L^{2}$ when the initial
 data are randomized.  Based upon Bourgain's high-low
frequency decomposition firstly used in \cite{CO} in probabilistic setting  and improved averaging
effects for the free evolution
of the randomized initial data, L$\ddot{u}$hrmann and  Mendelson \cite{LM} studied the  random data
 Cauchy theory for nonlinear wave equations of power-type on $\R^{3}.$
  Oh and his collaborators  \cite{BO-tms,BOP,OP,P} studied the probabilitistic Cauchy
  theory of the nonlinear Schr$\ddot{o}$dinger equations
 and wave equation in the supercritical regime. Recently, Nahmod and Staffilani \cite{NS} established
 the almost surely local well-posedness
 result of  the periodic 3D quintic nonlinear Schr\"odinger equation in the supercritical regime.

In this paper, motivated
 by  \cite{B1994,B-1996,BT2008-I,BT2008,BOP,OP, BO-tms}, we prove that the random data  Cauchy problem for (\ref{1.01})  is almost surely locally
well-posed in  $H^{s}(\R)$ with $s>\frac{17}{112}$, whose lower bound is below  the scaling critical
regularity $s_{crit}=\frac{3}{14}$ with the aid of  the Wiener randomization of the initial data and
improved local-in-time Strichartz estimate as well as some embedding Theorems.

We always assume that $b=\frac{1}{2}+\frac{\epsilon}{24}$ and
$c=\frac{1}{2}+\frac{\epsilon}{100}$. In this paper,
$C>0$ and $C^{'}>0$ denote two constants which may vary at different occurrence.
Let $\eta(t)$ be a smooth cutoff function
supported on $[-2,2]$, $\eta\equiv1$ on $[-1,1]$, and $\eta_{T}(t)=
\eta\left(\frac{t}{T}\right)$. For $x,\xi \in \R$, we denote by
$\mathscr{F}_{x}u(\xi)=\frac{1}{\sqrt{2\pi}} \int_{\SR}u(x)e^{-ix\xi}dx$
 the Fourier transformation with respect
to the space variable. We denote by
$\mathscr{F}_{t}u(\tau)=\frac{1}{\sqrt{2\pi}} \int_{\SR}u(t)e^{-it\tau}dt$
 the Fourier transformation with respect
to the time variable.
We denote by $\mathscr{F}u(\xi,\tau)=\frac{1}{2\pi} \int_{\SR^{2}}u(x,t)e^{-ix\xi-it\tau}dxdt$
 the Fourier transformation
with respect to the time and space variables. We define $\mathscr{F}_{x}J^{s}u(\xi)
=\langle \xi\rangle ^{s}\mathscr{F}_{x}u(\xi)$ and
$S(t)\phi=\frac{1}{\sqrt{2\pi}}\int_{\SR}e^{ix\xi}\mathscr{F}_{x}\phi(\xi)e^{it\xi^{3}}d\xi$.
Let $\R=\bigcup \limits_{n\in Z}Q_{n}$, where $Q_{n}=\left\{x\left|\right|x-n|\leq 1\right\}$.
Let $\supp \psi\subset [-1,1]$ and $\sum\limits_{n\in Z}\psi(\xi-n)=1$. We define
$\psi(D-n)u(x)=\int_{\SR}\psi(\xi-n)\mathscr{F}_{x}u(\xi)e^{i x\xi}d\xi$.
Obviously,
\begin{eqnarray*}
&&\sum\limits_{n\in Z}\psi(D-n)u(x)=\sum\limits_{n\in Z}\int_{\SR}\psi(\xi-n)
\mathscr{F}_{x}u(\xi)e^{i x\xi}d\xi\\&&=\int_{\SR}\sum_{n\in Z}\psi(\xi-n)
\mathscr{F}_{x}u(\xi)e^{i x\xi}d\xi=\int_{\SR}\mathscr{F}_{x}u(\xi)e^{i x\xi}d\xi=u(x).
\end{eqnarray*}
Let $0<p,q\leq \infty$ and $s\in \R$. The space $M_{s}^{p,q}$  consists of all
tempered distributions $u\in \mathscr{S}^{'}(\R)$
for which the (quasi) norm $\|u\|_{M_{s}^{p,q}(\SR)}=
\left\|\langle n\rangle ^{s}\|\psi(D-n)u\|_{L_{x}^{p}(\SR)}\right\|_{l_{n}^{q}(Z)}$.
 The  Sobolev space $H^{s}(\R)$ is defined as follows:
\begin{eqnarray*}
H^{s}(\R)=\left\{u\in \mathscr{S}^{'}(\R):\|u\|_{H^{s}(\SR)}=\|\langle \xi\rangle ^{s}\mathscr{F}_{x}u\|_{L_{\xi}^{2}}<\infty\right\}
\end{eqnarray*}
and
the   space $X_{s,b}(\R^{2})$ is defined as follows:
\begin{eqnarray*}
X_{s,b}(\R^{2})=\left\{u\in \mathscr{S}^{'}(\R^{2}):\|u\|_{X_{s,b}(\SR^{2})}=\|\langle \xi\rangle ^{s}\langle \tau-\xi^{3}\rangle^{b}\mathscr{F}u(\xi,\tau)\|_{L_{\xi}^{2}L_{\tau}^{2}}<\infty\right\}.
\end{eqnarray*}
These spaces were introduced
in the study of propagation of singularity in semilinear wave equation
by Rauch and Reed \cite{RR} and Beals \cite{B},
which  were used to systematically study nonlinear dispersive wave
problems by Bourgain \cite{Bourgain-GAFA93}. Moreover,
 Klainerman and Machedon \cite{KM} used similar ideas
in their study of the nonlinear wave equation.

Let $\left\{g_{n}\right\}_{n\in Z}$ be a sequence of
independent mean zero complex-valued
random variables on a probability space $(\Omega, \mathcal {F},P),$
where the real and imaginary part of $g_{n}$ are
independent and endowed with probability
distributions $\mu_{n}^{1}$ and $\mu_{n}^{2},$ respectively. For a function $\phi$ on $\R$,
we define the Wiener randomization of $\phi$  by
\begin{eqnarray}
\phi^{\omega}=\sum_{n\in Z}g_{n}(\omega)\psi(D-n)\phi.\label{1.02}
\end{eqnarray}
The Wiener randomization of $\phi$ does not improve the differentiability \cite{BT2008,BT2008-I} in Sobolev spaces,
however improves the integrability \cite{BOP}.
Thus, the advantage of
 the Wiener randomization is to make the critical  problem  in  Sobolev  space become subcritical in some sense.

The main result of this paper is as follows.
\begin{Theorem}\label{Thm1}
Let $s>\frac{17}{112}$, $\phi\in H^{s}(\R)$ and $\phi^{\omega}$ be its
randomization defined in (\ref{1.02}) and
\begin{eqnarray}
\left|\int_{\SR}e^{\gamma x}d\mu_{n}^{(j)}(x)\right|\leq e^{C\gamma^{2}}\label{1.04}
\end{eqnarray}
for all $\gamma\in \R,n\in Z,j=1,2.$ Then  (\ref{1.01})  is  almost
surely locally well-posed with respect
to the randomization
$\phi^{\omega}$ as initial data. More  precisely, there exist constants  $C,C^{'}>0$ and
$\sigma=\frac{3}{14}+2\epsilon$
such that for each $T\ll1$,
there exists an event $\Omega_{T}\subset \Omega$ with the following properties:

\noindent {\rm (i)}${\rm P}\left(\Omega\setminus\Omega_{T}\right)\leq C{\rm exp}
\left(-\frac{C^{'}}{T^{\frac{7\epsilon}{2400}}
\|\phi\|_{H^{s}}^{2}}\right).$

\noindent {\rm (ii)}For each $\omega\in \Omega_{T},$ there exists a unique
solution to (\ref{1.01})
with $u(x,0)=\phi^{\omega}$ in the class
$
S(t)\phi^{\omega}+C([-T,T]:H^{\sigma}(\R))\subset C([-T,T]:H^{s}(\R)).
$
\end{Theorem}

The rest of the paper is arranged as follows. In Section 2,  we give some
preliminaries. In Section 3, we show some multilinear estimates. In Section 4, we prove
the Theorem 1.1.

\bigskip

 \noindent{\large\bf 2. Preliminaries }

\setcounter{equation}{0}

\setcounter{Theorem}{0}

\setcounter{Lemma}{0}

\setcounter{section}{2}

In this section, we present some probabilistic lemmata and probabilistic Strichartz estimates,  linear estimates needed in this paper. We always denote that
$S(t)\phi=e^{-t\partial_{x}^{3}}=\frac{1}{2\sqrt{\pi}}\int_{\SR}e^{ix\xi}\mathscr{F}_{x}\phi(\xi)e^{it\xi^{3}}d\xi$.

\begin{Lemma}\label{Lemma2.1}
Let $\phi \in H^{s}(\R)$ and  $\phi ^{\omega}=\sum\limits_{n\in Z}g_{n}(\omega) \Psi(D-n)\phi$. Then, we have the probabilistic estimate
$
{\rm P} \left(\omega:\|\phi^{\omega}\|_{H^{s}}>\lambda\right)\leq Ce^{-C^{'}\lambda^{2}\|\phi\|_{H^{s}}^{-2}}
$
for all $\lambda>0$.
\end{Lemma}

For the proof of Lemma 2.1, we refer the readers to Lemma 2.2 of \cite{BOP}.

\begin{Lemma}\label{Lemma2.2}Let $\phi \in L^{2}(\R)$ and
$\phi ^{\omega}=\sum\limits_{n\in Z}g_{n}(\omega) \psi(D-n)\phi$.
Then,  for $2\leq q,r<\infty,$ there exist $C>0, C^{'}>0$ such that
$$
P\left(\omega:\left\|S(t)\phi^{\omega}\right\|_{L_{t}^{q}L_{x}^{r}([0,T]\times \SR)}>\lambda\right)
\leq C{\rm exp}\left(-C^{'}\frac{\lambda^{2}}{T^{2/q}\|\phi\|_{L^{2}}^{2}}\right)
$$
for all $T>0$ and $\lambda>0$.
\end{Lemma}
In particular, when $\lambda=T^{\theta}\|\phi\|_{L^{2}}$, we have
$
\left\|\langle\nabla\rangle^{s}S(t)\phi^{\omega}\right\|_{L_{t}^{q}L_{x}^{r}([0,T]\times \SR)}
\leq CT^{\theta}\|\phi\|_{H^{s}(\SR)}
$
outside a set of probability at most $C{\rm exp}\left(-C^{'}T^{2\theta-\frac{2}{q}}\right)$.

Lemma 2.2 can be proved  similarly to Proposition 1.3 of \cite{BOP}.

\begin{Lemma}\label{Lemma2.3}
Given $\epsilon>0$ and $T_{0}>0,$  there exist constants $C,C^{'}>0$, depending on $\epsilon$, such that
\begin{eqnarray}
P\left(\omega:\|S(t)\phi^{\omega}\|_{L_{t}^{\infty}([0,T_{0}];L_{x}^{\infty}(\SR))}>R\right)\leq C(1+T_{0}){\rm exp}\left(-C^{'}\frac{R^{2}}{\|\phi\|_{H^{\epsilon}(\SR)}^{2}}\right).\label{2.01}
\end{eqnarray}
\end{Lemma}
\noindent{\bf Proof.} By using a proof similar to (3.3) of  Lemma 3.4 of
 \cite{OP}, we have that
 \begin{eqnarray}
 P\left(\omega:\|S(t)\phi^{\omega}\|_{L_{t}^{\infty}([j,j+1];L_{x}^{\infty}(\SR))}>R\right)
 \leq C\left(-C^{'}\frac{R^{2}}{\|\phi\|_{H^{\epsilon}(\SR)}^{2}}\right)\label{2.02}
 \end{eqnarray}
for $j\in {\bf N}\cup 0$. Combining (\ref{2.02}) with the subadditivity, we have that
\begin{eqnarray}
&&P\left(\omega:\|S(t)\phi^{\omega}\|_{L_{t}^{\infty}([0,T_{0}];L_{x}^{\infty}(\SR))}>R\right)
\leq \sum_{j=0}^{[T_{0}]}P\left(\omega:\|S(t)\phi^{\omega}\|_{L_{t}^{\infty}([j,j+1];L_{x}^{\infty}(\SR))}>R\right)\nonumber\\
 &&\leq  C(1+T_{0}){\rm exp}\left(-C^{'}\frac{R^{2}}{\|\phi\|_{H^{\epsilon}(\SR)}^{2}}\right).\label{2.03}
\end{eqnarray}

We have completed the proof of Lemma 2.3.

\begin{Lemma}\label{Lemma2.4}
Let $T\in (0,1)$ and $b\in (\frac{1}{2},\frac{3}{2})$. Then, for $s\in \R$ and
$\theta \in [0,\frac{3}{2}-b)$, we have that
\begin{eqnarray}
&&\|\eta_{T}(t)S(t)\phi\|_{X_{s,b}(\SR^{2})}\leq CT^{\frac{1}{2}-b}\|\phi\|_{H^{s}(\SR)},\nonumber\\
&&\left\|\eta_{T}(t)\int_{0}^{t}S(t-\tau)F(\tau)d\tau\right\|_{X_{s,b}(\SR^{2})}
\leq CT^{\theta}\|F\|_{X_{s,b-1+\theta}(\SR^{2})}.\nonumber
\end{eqnarray}
\end{Lemma}

For the proof of Lemma 2.4, we refer the readers to \cite{Bourgain-GAFA93,KPVDuke,G}.

\begin{Lemma}\label{Lemma2.5}
Let $b_{1}>\frac{1}{2}\geq s\geq0,\tilde{b}>\frac{1}{6}+\frac{2s}{3}.$
 Then the following estimate holds true:
\begin{eqnarray*}
\left\|I^{s}I^{s}_{-}\left(u_{1},u_{2}\right)\right\|_{L_{xt}^{2}}\leq C
\|u_{1}\|_{X_{0,b_{1}}}\|u_{2}\|_{X_{0,\tilde{b}}}.
\end{eqnarray*}
\end{Lemma}

For the proof of Lemma 2.5, we refer the readers to \cite{G}.

\noindent In particular,
 we have that
\begin{eqnarray}
\left\|I^{\frac{1-\epsilon}{2}}I^{\frac{1-\epsilon}{2}}_{-}\left(u_{1},u_{2}\right)\right\|_{L_{xt}^{2}}\leq C\|u_{1}\|_{X_{0,\frac{1}{2}+\frac{\epsilon}{100}}}\|u_{2}\|_{X_{0,\frac{1}{2}-\frac{\epsilon}{12}}},\label{2.04}
\end{eqnarray}
\noindent
and
\begin{eqnarray}
\left\|I^{\frac{1}{2}}I^{\frac{1}{2}}_{-}\left(u_{1},u_{2}\right)\right\|_{L_{xt}^{2}}\leq C\|u_{1}\|_{X_{0,\frac{1}{2}+\frac{\epsilon}{100}}}\|u_{2}\|_{X_{0,\frac{1}{2}+\frac{\epsilon}{100}}}.\label{2.05}
\end{eqnarray}

\begin{Lemma}\label{Lemma2.6}
Let $0<\epsilon\ll1$ and $b=\frac{1}{2}+\frac{\epsilon}{24}$ and $3\leq l<7$. Then, we have that
\begin{eqnarray}
&&\|u\|_{L_{xt}^{\frac{8}{1+\epsilon}}}\leq C\|u\|_{X_{0,\frac{1}{2}-\frac{\epsilon}{12}}},\label{2.06}\\
&&\|u\|_{L_{xt}^{\frac{28}{2-7\epsilon}}}\leq C\|u\|_{X_{s_{c}+\epsilon,b}},\label{2.07}\\
&&\|u\|_{L_{xt}^{\frac{280}{17+7\epsilon}}}\leq C\|u\|_{X_{\frac{18-7\epsilon}{70},b}},\label{2.08}\\
&&\|u\|_{L_{xt}^{\frac{56(7-l)}{25-4l-7\epsilon(3-2l)}}}\leq C\|u\|_{X_{\frac{(8-l)(3+\epsilon)}{70},b}},\label{2.09}\\
&&\|u\|_{L_{xt}^{\frac{64}{7-\epsilon}}}\leq C\|u\|_{X_{\frac{1+\epsilon}{16},b}},\label{2.010}\\
&&\|u\|_{L_{xt}^{\frac{392}{45-84\epsilon}}}\leq C\|u\|_{X_{\frac{2+42\epsilon}{49},b}},\label{2.011}\\
&&\|u\|_{L_{xt}^{\frac{56}{4+77\epsilon}}}\leq C\|u\|_{X_{\frac{3-77\epsilon}{14},b}},\label{2.012}\\
&&\|u\|_{L_{xt}^{\frac{224}{13-70\epsilon}}}\leq C\|u\|_{X_{\frac{15+70\epsilon}{56},b}},\label{2.013}\\
&&\|u\|_{L_{xt}^{8}}\leq C\|u\|_{X_{0,b}},\label{2.014}\\
&&\|u\|_{L_{xt}^{\frac{32}{3-\epsilon}}}\leq C\|u\|_{X_{\frac{1+\epsilon}{8},b}},\label{2.015}\\
&&\|u\|_{L_{xt}^{\infty}}\leq C\|u\|_{X_{b,b}}.\label{2.016}
\end{eqnarray}
\end{Lemma}
\noindent {\bf Proof.}
It is easily checked that
\begin{eqnarray}
\|u\|_{L_{xt}^{8}}\leq C\|u\|_{X_{0,\frac{6-\epsilon}{12-4\epsilon}}}\label{2.017}
\end{eqnarray}
and
\begin{eqnarray}
\|u\|_{L_{xt}^{2}}\leq C\|u\|_{X_{0,0}}\label{2.018}.
\end{eqnarray}
Interpolating (\ref{2.017}) with (\ref{2.018}), we have that (\ref{2.06}) is valid.
By using the Sobolev embedding Theorem and  Theorem 2.4 of \cite{KPV1991},  we have that
\begin{eqnarray}
&&\hspace{-1.7cm}\|S(t)\phi\|_{L_{xt}^{\frac{28}{2-7\epsilon}}}\leq C
\|D_{x}^{\frac{3+14\epsilon}{56}}D_{t}^{\frac{3+14\epsilon}{56}}S(t)\phi\|_{L_{xt}^{8}}
\leq C\|D_{x}^{\frac{3+14\epsilon}{14}}S(t)\phi\|_{L_{xt}^{8}}\leq
C\|D_{x}^{\frac{3+14\epsilon}{14}}\phi\|_{L_{xt}^{2}}.\label{2.019}
\end{eqnarray}
Combining (\ref{2.019})   with a standard argument proof \cite{KPV1991}, we have
$
\|u\|_{L_{xt}^{\frac{28}{2-7\epsilon}}}\leq \|u\|_{X_{s_{c}+\epsilon,b}}.
$
Similarly, we have that (\ref{2.08})-(\ref{2.013}), (\ref{2.015}) are valid.
Obviously, (\ref{2.016}) is valid and
(\ref{2.014}) is in \cite{KPV1991}.

We have completed the proof of
Lemma 2.6.
\bigskip

\noindent{\large\bf 3. Multilinear estimates }

\setcounter{equation}{0}

 \setcounter{Theorem}{0}

\setcounter{Lemma}{0}

 \setcounter{section}{3}
 In this section, we present some crucial multilinear estimates which play an important role in establishing Theorem 1.1.

\begin{Lemma}\label{Lemma3.1}
Let $s=\frac{17}{112}+\epsilon$,  $\sigma=\frac{3}{14}+2\epsilon$ and $v_{j}=\eta v$ with $(1\leq j\leq 8,j\in N)$.
Then, we have that
\begin{eqnarray*}
\left|\int_{\SR}\int_{\SR}J^{\sigma}\partial_{x}\left(\prod_{j=1}^{8}v_{j}\right)h(x,t)dxdt\right|\leq C\left(\prod_{j=1}^{8}\|v_{j}\|_{X_{\sigma,b}}\right)\|h\|_{X_{0,\frac{1}{2}-\frac{\epsilon}{12}}}.
\end{eqnarray*}
\end{Lemma}
\noindent {\bf Proof.}
In this case, we divide the frequency  into
$$|\xi_{l}|\geq {\rm max}\left\{|\xi_{j}|,1\leq j\leq 8,j\neq l,j\in N\right\}(1\leq l\leq 8,l\in N).$$
Without loss of generality, we  assume that $|\xi_{1}|\geq |\xi_{2}|\geq\cdot\cdot\cdot\geq |\xi_{8}|.$

\noindent We define $$I_{1}=\left|\int_{\SR}\int_{\SR}J^{\sigma}\partial_{x}\left(\prod_{j=1}^{8}v_{j}\right)hdxdt\right|.$$
\noindent When $|\xi_{1}|\geq80 |\xi_{2}|,$   by using  the H\"older inequality, (\ref{2.04})-(\ref{2.05}) and
(\ref{2.016}),
we infer that
\begin{eqnarray*}
&&I_{1}\leq C
\left\|I^{1/2}I^{1/2}_{-}(J ^{\sigma}v_{1},J ^{\sigma}v_{2})\right\|_{L_{xt}^{2}}
\left(\prod_{j=4}^{8}\|J^{-\frac{2}{7}}v_{j}\|_{L_{xt}^{\infty}}\right)
\|I^{\frac{1-\epsilon}{2}}I^{\frac{1-\epsilon}{2}}_{-}(J ^{\sigma}v_{3},h)\|_{L_{xt}^{2}}\nonumber\\&&\leq
C\left(\prod_{j=1}^{8}\|v_{j}\|_{X_{\sigma,b}}\right)\|h\|_{X_{0,\frac{1}{2}-\frac{\epsilon}{12}}}.
\end{eqnarray*}
\noindent When $|\xi_{1}|\sim |\xi_{2}|\geq 80|\xi_{3}|,$ by using the  H\"older inequality, (\ref{2.05})-(\ref{2.08}),
we have that
\begin{eqnarray*}
&&I_{1}
\leq C\|v_{1}\|_{L_{xt}^{\frac{28}{2-7\epsilon}}}
\left\|I^{1/2}I^{1/2}_{-}(J^{\sigma}v_{2},J^{\sigma}v_{3})\right\|_{L_{xt}^{2}}
\left(\prod_{j=4}^{8}\|J^{-\frac{\sigma}{5}}v_{j}\|_{L_{xt}^{\frac{280}{17+7\epsilon}}}\right)
\|h\|_{L_{xt}^{\frac{8}{1+\epsilon}}}\nonumber\\&&\leq
C\left(\prod_{j=1}^{8}\|v_{j}\|_{X_{\sigma,b}}\right)\|h\|_{X_{0,\frac{1}{2}-\frac{\epsilon}{12}}}.
\end{eqnarray*}
\noindent When $|\xi_{1}|\sim |\xi_{l}|\geq 80|\xi_{l+1}|(3\leq l<7,l\in N),$  by using the  H\"older inequality together with
(\ref{2.05})-(\ref{2.07}) and (\ref{2.09}), we have that
\begin{eqnarray*}
&&I_{1}\leq C
\left\|I^{1/2}I^{1/2}_{-}(J^{\sigma}v_{1},J^{\sigma}v_{l+1})\right\|_{L_{xt}^{2}}
\left(\prod_{j=2}^{l}\|v_{j}\|_{L_{xt}^{\frac{28}{2-7\epsilon}}}\right)\nonumber\\&&\qquad\times
\left(\prod_{j=l+2}^{8}\|J^{-\frac{\sigma}{7-l}}v_{j}\|_{L_{xt}^{\frac{56(7-l)}{25-4l-7\epsilon(3-2l)}}}\right)
\|h\|_{L_{xt}^{\frac{8}{1+\epsilon}}}\leq C\left(\prod_{j=1}^{8}
\|v_{j}\|_{X_{\sigma,b}}\right)\|h\|_{X_{0,\frac{1}{2}-\frac{\epsilon}{12}}}.
\end{eqnarray*}
\noindent When $|\xi_{1}|\sim |\xi_{7}|\geq 80|\xi_{8}|,$  by using  the  H\"older inequality, (\ref{2.06}),  (\ref{2.011}),
(\ref{2.012}), we obtain that
\begin{eqnarray*}
&&I_{1}\leq C
\left(\prod_{j=1}^{7}\|J^{\frac{17+28\epsilon}{98}}v_{j}\|_{L_{xt}^{\frac{392}{45-84\epsilon}}}\right)
\|v_{8}\|_{L_{xt}^{\frac{56}{4+77\epsilon}}}
\|h\|_{L_{xt}^{\frac{8}{1+\epsilon}}}\leq C\left(\prod_{j=1}^{8}\|v_{j}\|_{X_{\sigma,b}}\right)
\|h\|_{X_{0,\frac{1}{2}-\frac{\epsilon}{12}}}.
\end{eqnarray*}
When $|\xi_{1}|\sim|\xi_{8}|$, by using the  H\"older inequality, (\ref{2.06}),  (\ref{2.010}), we have that
\begin{eqnarray*}
&&I_{1}\leq C
\left(\prod_{j=1}^{8}\|J^{\sigma-\frac{1+\epsilon}{16}}v_{j}\|_{L_{xt}^{\frac{64}{7-\epsilon}}}\right)
\|h\|_{L_{xt}^{\frac{8}{1+\epsilon}}}\leq  C\left(\prod_{j=1}^{8}\|v_{j}\|_{X_{\sigma,b}}\right)
\|h\|_{X_{0,\frac{1}{2}-\frac{\epsilon}{12}}}.
\end{eqnarray*}

This completes the proof of Lemma 3.1.

\begin{Lemma}\label{Lemma3.2}
Let
  $s\geq\frac{17}{112}+\epsilon$ and $\sigma=\frac{3}{14}+2\epsilon$ and
   $z_{j}(t)=\eta_{T}(t)S(t)\phi^{\omega}(1\leq j\leq8,j\in N)$. Then, we have that
\begin{eqnarray*}
\left|\int_{\SR}\int_{\SR}J^{\sigma}\partial_{x}\left(\prod_{j=1}^{8}z_{j}\right)hdxdt\right|\leq CT^{-\frac{3\epsilon}{100}}R^{8}\|h\|_{X_{0,\frac{1}{2}-\frac{\epsilon}{12}}},
\end{eqnarray*}
outside a set of probability at most
$
 C{\rm exp}\left(-C^{'}\frac{R^{2}}{\|\phi\|_{H^{s}}^{2}}\right).
$
\end{Lemma}
\noindent{\bf Proof.}We dyadically decompose $z_{j}$ with $(1\leq j\leq8,j\in N)$ and
$h$ such that frequency supports are
$\left\{|\xi_{j}|\sim N_{j}\right\}$ for some dyadically $N_{j}\geq 1$ and we still
denote them by $z_{j}$ with $(1\leq j\leq 8,j\in N)$ and $h$.
In this case, we divide the frequency  into
$$|\xi_{l}|\geq {\rm max}\left\{|\xi_{j}|,1\leq j\leq 8,j\neq l,j\in Z\right\}(1\leq l\leq 8,l\in N).$$
Without loss of generality, we  assume that $|\xi_{1}|\geq |\xi_{2}|\geq\cdot\cdot\cdot\geq |\xi_{8}|.$

\noindent We define $$I_{2}=\left|\int_{\SR}\int_{\SR}J^{\sigma}\partial_{x}\left(\prod_{j=1}^{8}z_{j}\right)hdxdt\right|.$$
\noindent When $|\xi_{1}|\geq80|\xi_{2}|,$ since   $\sigma-s-1+\epsilon<0$,  by using the  H\"older inequality,
(\ref{2.04})-(\ref{2.05}), Lemmas 2.3, 2.4, 2.1,
we have that
\begin{eqnarray*}
&&I_{2}\leq CN_{1}^{\sigma-s-1+\epsilon}
\left\|I^{1/2}I^{1/2}_{-}(J ^{s}z_{1},J ^{s}z_{2})\right\|_{L_{xt}^{2}}
\left(\prod_{j=4}^{8}\|z_{j}\|_{L_{xt}^{\infty}}\right)
\|I^{\frac{1-\epsilon}{2}}I^{\frac{1-\epsilon}{2}}_{-}(J ^{s}z_{3},h)\|_{L_{xt}^{2}}\nonumber\\&&\leq
CN_{1}^{\sigma-s-1+\epsilon}\left(\prod_{j=1}^{3}\|z_{j}\|_{X_{s,c}}\right)
\left(\prod_{j=4}^{8}\|z_{j}\|_{L_{xt}^{\infty}}\right)
\|h\|_{X_{0,\frac{1}{2}-\frac{\epsilon}{12}}}\nonumber\\&&\leq
CN_{1}^{\sigma-s-1+\epsilon}T^{-\frac{3\epsilon}{100}}\left(\prod_{j=1}^{3}\|P_{N_{j}}
\phi^{\omega}\|_{H^{s}}\right)\left(\prod_{j=4}^{8}\|z_{j}\|_{L_{xt}^{\infty}}\right)
\|h\|_{X_{0,\frac{1}{2}-\frac{\epsilon}{12}}}\nonumber\\&&\leq C
T^{-\frac{3\epsilon}{100}}R^{8}\|h\|_{X_{0,\frac{1}{2}-\frac{\epsilon}{12}}},
\end{eqnarray*}
outside a set of probability at most
$
C{\rm exp}\left(-C^{'}\frac{R^{2}}{\|\phi\|_{H^{s}}^{2}}\right).
$

\noindent When $|\xi_{1}|\sim |\xi_{2}|\geq 80|\xi_{3}|,$ by using the
 H\"older inequality and $-2s+\sigma <0$, (\ref{2.05})-(\ref{2.06}), Lemmas 2.4, 2.1, 2.2,
we have that
\begin{eqnarray*}
&&I_{2}\leq CN_{1}^{-2s+\sigma}\|J ^{s}z_{1}\|_{L_{xt}^{\frac{48}{3-\epsilon}}}
\left\|I^{1/2}I^{1/2}_{-}(J ^{s}z_{2},z_{3})\right\|_{L_{xt}^{2}}
\left(\prod_{j=4}^{8}\|z_{j}\|_{L_{xt}^{\frac{48}{3-\epsilon}}}\right)
\|h\|_{L_{xt}^{\frac{8}{1+\epsilon}}}\nonumber\\&&\leq CN_{1}^{-2s+\sigma}
\|J ^{s}z_{1}\|_{L_{xt}^{\frac{48}{3-\epsilon}}}\left(\prod_{j=2}^{3}
\|z_{j}\|_{X_{s,c}}\right)
\left(\prod_{j=4}^{8}\|z_{j}\|_{L_{xt}^{\frac{48}{3-\epsilon}}}\right)
\|h\|_{X_{0,\frac{1}{2}-\frac{\epsilon}{12}}}\nonumber\\&&\leq
CN_{1}^{-2s+\sigma}T^{-\frac{\epsilon}{50}}\|J ^{s}z_{1}\|_{L_{xt}^{\frac{48}{3-\epsilon}}}
\left(\prod_{j=2}^{3}\|P_{N_{j}}\phi^{\omega}\|_{H^{s}}\right)
\left(\prod_{j=4}^{8}\|z_{j}\|_{L_{xt}^{\frac{48}{3-\epsilon}}}\right)
\|h\|_{X_{0,\frac{1}{2}-\frac{\epsilon}{12}}}\nonumber\\&&\leq
CT^{-\frac{\epsilon}{50}}R^{8}\|h\|_{X_{0,\frac{1}{2}-\frac{\epsilon}{12}}},
\end{eqnarray*}
outside a set of probability at most
$
 C{\rm exp}\left(-C^{'}\frac{R^{2}}{\|\phi\|_{H^{s}}^{2}}\right).
$

\noindent When $|\xi_{1}|\sim |\xi_{l}|\geq 80|\xi_{l+1}|(3\leq l\leq7,l\in N),$
since $-ls+\sigma<0$, by using the H\"older inequality,
(\ref{2.05}),  (\ref{2.06}) and Lemmas 2.4, 2.1, 2.2,   we have that
\begin{eqnarray*}
&&I_{2}\leq CN_{1}^{-ls+\sigma}
\left\|I^{1/2}I^{1/2}_{-}(J ^{s}z_{1},z_{l+1})\right\|_{L_{xt}^{2}}
\left(\prod_{j=2}^{l}\|J^{s}z_{j}\|_{L_{xt}^{\frac{48}{3-\epsilon}}}\right)
\left(\prod_{j=l+2}^{8}\|z_{j}\|_{L_{xt}^{\frac{48}{3-\epsilon}}}\right)
\|h\|_{L_{xt}^{\frac{8}{1+\epsilon}}}\nonumber\\&&\leq CN_{1}^{-ls+\sigma}
\|z_{1}\|_{X_{s,c}}
\|z_{l+1}\|_{X_{0,c}}
\left(\prod_{j=2}^{l}\|J^{s}z_{j}\|_{L_{xt}^{\frac{48}{3-\epsilon}}}\right)
\left(\prod_{j=l+2}^{8}\|z_{j}\|_{L_{xt}^{\frac{48}{3-\epsilon}}}\right)
\|h\|_{X_{0,\frac{1}{2}-\frac{\epsilon}{12}}}\nonumber\\&&\leq C
N_{1}^{-ls+\sigma}\|z_{1}\|_{X_{s,c}}
\|z_{l+1}\|_{X_{0,c}}
\left(\prod_{j=2}^{l}\|J^{s}z_{j}\|_{L_{xt}^{\frac{48}{3-\epsilon}}}\right)
\left(\prod_{j=l+2}^{8}\|z_{j}\|_{L_{xt}^{\frac{48}{3-\epsilon}}}\right)
\|h\|_{X_{0,\frac{1}{2}-\frac{\epsilon}{12}}}\nonumber\\&&\leq
CN_{1}^{-ls+\sigma}T^{-\frac{\epsilon}{50}}\|P_{N_{1}}\phi^{\omega}\|_{H^{s}}
\left(\prod_{j=2}^{l}\|J^{s}z_{j}\|_{L_{xt}^{\frac{48}{3-\epsilon}}}\right)
\|P_{N_{l+1}}\phi^{\omega}\|_{L^{2}}
\left(\prod_{j=l+2}^{8}\|z_{j}\|_{L_{xt}^{\frac{48}{3-\epsilon}}}\right)\nonumber\\&&\qquad\times
\|h\|_{X_{0,\frac{1}{2}-\frac{\epsilon}{12}}}\leq
CT^{-\frac{\epsilon}{50}}R^{8}\|h\|_{X_{0,\frac{1}{2}-\frac{\epsilon}{12}}},
\end{eqnarray*}
outside a set of probability  at most
$
C{\rm exp}\left(-C^{'}\frac{R^{2}}{\|\phi\|_{H^{s}}^{2}}\right).
$

\noindent When $|\xi_{1}|\sim |\xi_{8}|$, since  $-8s+1+\sigma<0,$
by using the H\"older  inequality,
 (\ref{2.06}) and Lemma 2.2, we get that
\begin{eqnarray*}
&&I_{2}\leq CN_{1}^{-8s+1+\sigma}
\left(\prod_{j=1}^{8}\|J^{s}z_{j}\|_{L_{xt}^{\frac{64}{7-\epsilon}}}\right)
\|h\|_{L_{xt}^{\frac{8}{1+\epsilon}}}\nonumber\\&&\leq CN_{1}^{-8s+1+\sigma}
\left(\prod_{j=1}^{8}\|J^{s}z_{j}\|_{L_{xt}^{\frac{64}{7-\epsilon}}}\right)
\|h\|_{X_{0,\frac{1}{2}-\frac{\epsilon}{12}}}\leq
CR^{8}\|h\|_{X_{0,\frac{1}{2}-\frac{\epsilon}{12}}},
\end{eqnarray*}
outside a set of probability at most
$
C{\rm exp}\left(-C^{'}\frac{R^{2}}{T^{\frac{7-\epsilon}{32}}\|\phi\|_{H^{s}}^{2}}\right).
$

We have completed the proof of Lemma 3.2.

\noindent {\bf Remark 2:} The last case of Lemma 3.2 requires that $s>
\frac{1}{8}\left(1+\sigma\right)=\frac{17+14\epsilon}{112}.$
\begin{Lemma}\label{Lemma3.3}
Let $s\geq \frac{17}{112}+\epsilon$ and $\sigma=\frac{3}{14}+2\epsilon$, $v_{j}=\eta (t) v$
with $(1\leq j\leq 7,j\in N)$ and $z_{8}=\eta_{T}(t)S(t)\phi^{\omega}$. Then, we have that
\begin{eqnarray*}
\left|\int_{\SR}\int_{\SR}J^{\sigma}\partial_{x}
\left(\left(\prod_{j=1}^{7}v_{j}\right)z_{8}\right)h(x,t)dxdt\right|\leq CRT^{-\frac{\epsilon}{100}}
\left(\prod_{j=1}^{7}\|v_{j}\|_{X_{\sigma,b}}\right)
\|h\|_{X_{0,\frac{1}{2}-\frac{\epsilon}{12}}},
\end{eqnarray*}
outside a set of probability at most
$
C{\rm exp}\left(-C^{'}\frac{R^{2}}{\|\phi\|_{H^{s}}^{2}}\right).
$
\end{Lemma}
\noindent {\bf Proof.} We dyadically decompose $v_{j}$ with $(1\leq j\leq7,j\in N)$,
$z_{8}$ and $h$ such that frequency supports are
$\left\{|\xi_{j}|\sim N_{j}\right\}$ for some dyadically $N_{j}\geq 1$ and we still
denote them by $v_{j}$ with $(1\leq j\leq 7,j\in N)$, $z_{8}$ and $h$.
In this case, we divide the frequency  into
$$|\xi_{l}|\geq {\rm max}\left\{|\xi_{j}|,1\leq j\leq 7,j\neq l,j\in N\right\}(1\leq l\leq 7,l\in N).$$
Without loss of generality, we can assume that $|\xi_{1}|\geq |\xi_{2}|\geq\cdot\cdot\cdot\geq |\xi_{7}|.$

\noindent We define $$I_{3}=\left|\int_{\SR}\int_{\SR}J^{\sigma}\partial_{x}
\left(\left(\prod_{j=1}^{7}v_{j}\right)z_{8}\right)h(x,t)dxdt\right|.$$
\noindent When $|\xi_{1}|\geq80 |\xi_{2}|,$ $|\xi_{1}|\geq 80|\xi_{8}|$,  by using the
H\"older inequality,   Lemma 2.3 and  (\ref{2.04}), (\ref{2.05}) and (\ref{2.016}),
we have that
\begin{eqnarray*}
&&\hspace{-1cm}I_{3}\leq CN_{1}^{-\frac{1}{10}}
\left\|I^{1/2}I^{1/2}_{-}(J ^{\sigma}v_{1},J ^{\sigma}v_{2})\right\|_{L_{xt}^{2}}
\left\|I^{\frac{1-\epsilon}{2}}I^{\frac{1-\epsilon}{2}}_{-}(J ^{\sigma}v_{3},h)\right\|_{L_{xt}^{2}}
\left(\prod_{j=4}^{7}\|J^{-\frac{2}{7}}v_{j}\|_{L_{xt}^{\infty}}\right)
\|z_{8}\|_{L_{xt}^{\infty}}
\nonumber\\&&\leq
CN_{1}^{-\frac{1}{10}}\left(\prod_{j=1}^{7}\|v_{j}\|_{X_{\sigma,b}}\right)\|z_{8}\|_{L_{xt}^{\infty}}
\|h\|_{X_{0,\frac{1}{2}-\frac{\epsilon}{12}}}\leq CR\left(\prod_{j=1}^{7}\|v_{j}\|_{X_{\sigma,b}}\right)
\|h\|_{X_{0,\frac{1}{2}-\frac{\epsilon}{12}}},
\end{eqnarray*}
outside a set of probability at most
$
C{\rm exp}\left(-C^{'}\frac{R^{2}}{\|\phi\|_{H^{\epsilon}}^{2}}\right).
$

\noindent When $|\xi_{1}|\geq80 |\xi_{2}|,$ $\frac{|\xi_{8}|}{80}\leq|\xi_{1}|\leq 80|\xi_{8}|$,
by using the H\"older inequality, (\ref{2.05})-(\ref{2.06}),  (\ref{2.08})
and Lemma 2.2,
we have that
\begin{eqnarray*}
&&I_{3}\leq C
\left\|I^{1/2}I^{1/2}_{-}(J ^{\sigma}v_{1},J ^{\sigma}v_{2})\right\|_{L_{xt}^{2}}\left(\prod_{j=3}^{7}
\|J^{-\frac{\sigma}{5}}v_{j}\|_{L_{xt}^{\frac{280}{17+7\epsilon}}}\right)
\|z_{8}\|_{L_{xt}^{\frac{28}{2-7\epsilon}}}
\|h\|_{L_{xt}^{\frac{8}{1+\epsilon}}}\nonumber\\&&\leq
C\left(\prod_{j=1}^{7}\|v_{j}\|_{X_{\sigma,b}}\right)
\|z_{8}\|_{L_{xt}^{\frac{28}{2-7\epsilon}}}\|h\|_{X_{0,\frac{1}{2}-\frac{\epsilon}{12}}}
\leq CR\left(\prod_{j=1}^{7}\|v_{j}\|_{X_{\sigma,b}}\right)
\|h\|_{X_{0,\frac{1}{2}-\frac{\epsilon}{12}}},
\end{eqnarray*}
outside a set of probability at most
$
C{\rm exp}\left(-C^{'}\frac{R^{2}}{T^{\frac{2-7\epsilon}{14}}\|\phi\|_{H^{\epsilon}}^{2}}\right).
$

\noindent
When $\frac{|\xi_{8}|}{80}\geq|\xi_{1}|\geq80 |\xi_{2}|,$ $|\xi_{1}|\leq|\xi_{8}|^{\frac{1}{3}},$
by using the H\"older inequality, (\ref{2.04})-(\ref{2.05}) and  (\ref{2.016}), Lemmas 2.4, 2.1,
we have that
\begin{eqnarray*}
&&I_{3}\leq CN_{8}^{-\frac{1}{2}}
\left\|I^{1/2}I^{1/2}_{-}(J^{s}z_{8},J ^{\sigma}v_{1})\right\|_{L_{xt}^{2}}
\left(\prod_{j=3}^{7}\|J^{-\frac{2}{7}}v_{j}\|_{L_{xt}^{\infty}}\right)
\|I^{\frac{1-\epsilon}{2}}I^{\frac{1-\epsilon}{2}}_{-}(J ^{\sigma}v_{2},h)\|_{L_{xt}^{2}}\nonumber\\&&\leq
CN_{8}^{-\frac{1}{2}}\left(\prod_{j=1}^{7}\|v_{j}\|_{X_{\sigma,b}}\right)\|z_{8}\|_{X_{s,c}}
\|h\|_{X_{0,\frac{1}{2}-\frac{\epsilon}{12}}}\nonumber\\
&&\leq CT^{-\frac{\epsilon}{100}}N_{8}^{-\frac{1}{2}}\left(\prod_{j=1}^{7}
\|v_{j}\|_{X_{\sigma,b}}\right)\|P_{N_{8}}\phi^{\omega}\|_{H^{s}}
\|h\|_{X_{0,\frac{1}{2}-\frac{\epsilon}{12}}}\leq
CT^{-\frac{\epsilon}{100}}R\left(\prod_{j=1}^{7}
\|v_{j}\|_{X_{\sigma,b}}\right),
\end{eqnarray*}
outside a set of probability at most
$
C{\rm exp}\left(-C^{'}\frac{R^{2}}{\|\phi\|_{H^{s}}^{2}}\right).
$

\noindent
When $\frac{|\xi_{8}|}{80}\geq|\xi_{1}|\geq80 |\xi_{2}|,$ $|\xi_{1}|\geq|\xi_{8}|^{\frac{1}{3}},$
by using the H\"older inequality, (\ref{2.04})-(\ref{2.05}),
Lemma 2.3, and (\ref{2.016}),
we have that
\begin{eqnarray*}
&&I_{3}\leq CN_{8}^{-\frac{1}{10}}
\left\|I^{1/2}I^{1/2}_{-}(J ^{\sigma}v_{1},J ^{\sigma}v_{2})\right\|_{L_{xt}^{2}}
\|I^{\frac{1-\epsilon}{2}}I^{\frac{1-\epsilon}{2}}_{-}(J ^{\sigma}v_{3},h)\|_{L_{xt}^{2}}\nonumber\\&&\qquad\left(\prod_{j=4}^{7}
\|J^{-\frac{2}{7}}v_{j}\|_{L_{xt}^{\infty}}\right)
\|J^{s}z_{8}\|_{L_{xt}^{\infty}}
\nonumber\\&&\leq
CN_{8}^{-\frac{1}{10}}\left(\prod_{j=1}^{7}\|v_{j}\|_{X_{\sigma,b}}\right)\|J^{s}z_{8}\|_{L_{xt}^{\infty}}
\|h\|_{X_{0,\frac{1}{2}-\frac{\epsilon}{12}}}\leq CR\left(\prod_{j=1}^{7}\|v_{j}\|_{X_{\sigma,b}}\right)
\|h\|_{X_{0,\frac{1}{2}-\frac{\epsilon}{12}}},
\end{eqnarray*}
outside a set of probability at most
$
C{\rm exp}\left(-C^{'}\frac{R^{2}}{\|\phi\|_{H^{s}}^{2}}\right).
$

\noindent When $|\xi_{1}|\sim |\xi_{2}|\geq 80|\xi_{3}|,$ $|\xi_{1}|\geq \frac{|\xi_{8}|}{80},$
by using the H\"older inequality,   (\ref{2.05})-(\ref{2.07}) and (\ref{2.013}), Lemma 2.2,
we have that
\begin{eqnarray*}
&&I_{3}\leq C\|v_{1}\|_{L_{xt}^{\frac{28}{2-7\epsilon}}}
\left\|I^{1/2}I^{1/2}_{-}(J^{\sigma}v_{2},J^{\sigma}v_{3})\right\|_{L_{xt}^{2}}
\left(\prod_{j=4}^{7}\|J^{-\frac{\sigma}{4}}v_{j}\|_{L_{xt}^{\frac{224}{13-70\epsilon}}}\right)
\|z_{8}\|_{L_{xt}^{\frac{56}{4+77\epsilon}}}
\|h\|_{L_{xt}^{\frac{8}{1+\epsilon}}}\nonumber\\&&\leq
C\left(\prod_{j=1}^{7}\|v_{j}\|_{X_{\sigma,b}}\right)\|z_{8}\|_{L_{xt}^{\frac{56}{4+77\epsilon}}}
\|h\|_{X_{0,\frac{1}{2}-\frac{\epsilon}{12}}}\leq CR\left(\prod_{j=1}^{7}\|v_{j}\|_{X_{\sigma,b}}\right)
\|h\|_{X_{0,\frac{1}{2}-\frac{\epsilon}{12}}},
\end{eqnarray*}
outside a set of probability at most
$
C{\rm exp}\left(-C^{'}\frac{R^{2}}{T^{\frac{4+77\epsilon}{28}}\|\phi\|_{H^{\epsilon}}^{2}}\right).
$

\noindent When $\frac{|\xi_{8}|}{80}\geq|\xi_{1}|\sim |\xi_{2}|\geq 80|\xi_{3}|,$
$|\xi_{1}|\leq|\xi_{8}|^{\frac{1}{3}},$  this case can be proved  similarly to case
 $\frac{|\xi_{8}|}{80}\geq|\xi_{1}|\geq 80|\xi_{2}|,$
$|\xi_{1}|\leq|\xi_{8}|^{\frac{1}{3}}$ of  Lemma 3.3.

\noindent When $\frac{|\xi_{8}|}{80}\geq|\xi_{1}|\sim |\xi_{2}|\geq 80|\xi_{3}|,$
$|\xi_{1}|\geq|\xi_{8}|^{\frac{1}{3}},$  this case  can  be  proved similarly to case
$\frac{|\xi_{8}|}{80}\geq|\xi_{1}|\geq 80|\xi_{2}|,$
$|\xi_{1}|\geq|\xi_{8}|^{\frac{1}{3}}$ of  Lemma 3.3.

\noindent When $|\xi_{1}|\sim |\xi_{l}|\geq 80|\xi_{l+1}|(3\leq l<6,l\in N),$
 $|\xi_{1}|\geq\frac{|\xi_{8}|}{80},$
this case can be proved  similarly to  case $|\xi_{1}|\sim |\xi_{2}|\geq 80|\xi_{3}|,$
 $|\xi_{1}|\geq\frac{|\xi_{8}|}{80},$ of Lemma 3.3.

\noindent When $\frac{|\xi_{8}|}{80}\geq|\xi_{1}|\sim |\xi_{l}|\geq 80|\xi_{l+1}|(3\leq l<6,l\in N),$
$|\xi_{1}|\leq|\xi_{8}|^{\frac{1}{3}},$
 this case can be proved  similarly to  case $\frac{|\xi_{8}|}{80}\geq|\xi_{1}|\geq 80|\xi_{2}|,$
$|\xi_{1}|\leq|\xi_{8}|^{\frac{1}{3}}$ of  Lemma 3.3.

\noindent When $\frac{|\xi_{8}|}{80}\geq|\xi_{1}|\sim |\xi_{l}|\geq 80|\xi_{l+1}|(3\leq l<6,l\in N),$
$|\xi_{1}|\geq|\xi_{8}|^{\frac{1}{3}},$
this case can be proved  similarly to  to case $\frac{|\xi_{8}|}{80}\geq|\xi_{1}|\geq 80|\xi_{2}|,$
$|\xi_{1}|\geq|\xi_{8}|^{\frac{1}{3}}$ of  Lemma 3.3.

\noindent When $|\xi_{1}|\sim |\xi_{6}|\geq 80|\xi_{7}|,$ $|\xi_{1}|\geq\frac{|\xi_{8}|}{80},$
by using the H\"older inequality, (\ref{2.05})-(\ref{2.07}),
   (\ref{2.016}) and Lemma 2.2, we obtain that
\begin{eqnarray*}
&&I_{3}
\leq C\left(\prod_{j=1}^{5}\|v_{j}\|_{L_{xt}^{\frac{28}{2-7\epsilon}}}\right)
\left\|I^{1/2}I^{1/2}_{-}(J^{\sigma}v_{6},J^{\sigma}v_{7})\right\|_{L_{xt}^{2}}
\|z_{8}\|_{L_{xt}^{\frac{56}{1+63\epsilon}}}
\|h\|_{L_{xt}^{\frac{8}{1+\epsilon}}}\nonumber\\&&\leq C
\left(\prod_{j=1}^{7}\|v_{j}\|_{X_{\sigma,b}}\right)
\|z_{8}\|_{L_{xt}^{\frac{56}{1+63\epsilon}}}\|h\|_{X_{0,\frac{1}{2}-\frac{\epsilon}{12}}}\leq
CR\left(\prod_{j=1}^{7}\|v_{j}\|_{X_{\sigma,b}}\right)
\|h\|_{X_{0,\frac{1}{2}-\frac{\epsilon}{12}}},
\end{eqnarray*}
outside a set of probability at most
$
C{\rm exp}\left(-C^{'}\frac{R^{2}}{T^{\frac{1+63\epsilon}{28}}\|\phi\|_{H^{\epsilon}}^{2}}\right).
$

\noindent When $\frac{|\xi_{8}|}{80}\geq|\xi_{1}|\sim |\xi_{6}|\geq 80|\xi_{7}|,$
$|\xi_{1}|\leq|\xi_{8}|^{\frac{1}{3}},$
this case can be proved  similarly to   case $\frac{|\xi_{8}|}{80}\geq|\xi_{1}|\geq 80|\xi_{2}|,$
$|\xi_{1}|\leq|\xi_{8}|^{\frac{1}{3}}$ of  Lemma 3.3.

\noindent When $\frac{|\xi_{8}|}{80}\geq|\xi_{1}|\sim |\xi_{6}|\geq 80|\xi_{7}|,$
$|\xi_{1}|\geq|\xi_{8}|^{\frac{1}{3}},$ this case can be proved similarly to case
$\frac{|\xi_{8}|}{80}\geq|\xi_{1}|\geq 80|\xi_{2}|,$
$|\xi_{1}|\geq|\xi_{8}|^{\frac{1}{3}}$ of  Lemma 3.3.

\noindent When $|\xi_{1}|\sim |\xi_{7}|,$  $|\xi_{1}|\geq \frac{|\xi_{8}|}{80},$
 by using the H\"older inequality, (\ref{2.06}),  (\ref{2.09}) and Lemma 2.2,  we conclude that
\begin{eqnarray*}
&&I_{3}\leq CN_{1}^{-6\epsilon}
\left(\prod_{j=1}^{7}\|J^{\sigma-\frac{2+42\epsilon}{49}}v_{j}\|_{L_{xt}^{\frac{392}{45-84\epsilon}}}\right)
\|z_{8}\|_{L_{xt}^{\frac{56}{4+77\epsilon}}}\|h\|_{L_{xt}^{\frac{8}{1+\epsilon}}}
\nonumber\\&&\leq CN_{1}^{-6\epsilon}
\left(\prod_{j=1}^{7}\|v_{j}\|_{X_{\sigma,b}}\right)
\|z_{8}\|_{L_{xt}^{\frac{56}{4+77\epsilon}}}\|h\|_{X_{0,\frac{1}{2}-\frac{\epsilon}{12}}}\leq CR
\left(\prod_{j=1}^{7}\|v_{j}\|_{X_{\sigma,b}}\right)
\|h\|_{X_{0,\frac{1}{2}-\frac{\epsilon}{12}}},
\end{eqnarray*}
outside a set of probability at most
$
C{\rm exp}\left(-C^{'}\frac{R^{2}}{T^{\frac{4+77\epsilon}{28}}\|\phi\|_{H^{\epsilon}}^{2}}\right).
$

\noindent When $\frac{|\xi_{8}|}{80}\geq|\xi_{1}|\sim |\xi_{7}|,$  $|\xi_{1}|\leq|\xi_{8}|^{\frac{1}{3}},$
this case can be proved similarly to case $\frac{|\xi_{8}|}{80}\geq|\xi_{1}|\geq 80|\xi_{2}|,$
$|\xi_{1}|\leq|\xi_{8}|^{\frac{1}{3}}$ of  Lemma 3.3.

\noindent When $\frac{|\xi_{8}|}{80}\geq|\xi_{1}|\sim |\xi_{7}|,$  $|\xi_{1}|\geq|\xi_{8}|^{\frac{1}{3}},$
this case can be proved similarly to case $\frac{|\xi_{8}|}{80}\geq|\xi_{1}|\geq 80|\xi_{2}|,$
$|\xi_{1}|\geq|\xi_{8}|^{\frac{1}{3}}$ of  Lemma 3.3.

\noindent
When $|\xi_{1}|\sim |\xi_{8}|,$ this case can be proved similarly to case
$|\xi_{1}|\sim |\xi_{7}|,$  $|\xi_{1}|\geq \frac{|\xi_{8}|}{80}$ of Lemma 3.3.

This completes the proof of Lemma 3.3.

\begin{Lemma}\label{Lemma3.4}
Let $s\geq\frac{17}{112}+\epsilon$ and $\sigma=\frac{3}{14}+2\epsilon$ and $v_{j}=\eta v(1\leq j\leq6,j\in N)$,
 $z_{j}=\eta_{T}S(t)\phi^{\omega}(j=7,8)$.
 Then, we have that
\begin{eqnarray*}
\left|\int_{\SR}\int_{\SR}J^{\sigma}\partial_{x}\left(\left(\prod_{j=1}^{6}v_{j}\right)z_{7}z_{8}\right)h(x,t)dxdt\right|\leq CT^{-\frac{\epsilon}{100}}R^{2}\left(\prod_{j=1}^{6}\|v_{j}\|_{X_{\sigma,b}}\right)\|h\|_{X_{0,\frac{1}{2}-\frac{\epsilon}{12}}},
\end{eqnarray*}
outside a set of probability  at most
$
C{\rm exp}\left(-C^{'}\frac{R^{2}}{\|\phi\|_{H^{s}}^{2}}\right).
$
\end{Lemma}
\noindent {\bf Proof.} We dyadically decompose $z_{j}$ with $( j=7,8)$ and $v_{j}$
with $(1\leq j\leq 6,j\in N)$ and $h$ such that frequency supports are
$\left\{|\xi_{j}|\sim N_{j}\right\}$ for some dyadically $N_{j}\geq 1$ and we still denote them
by $v_{j}$ with $(1\leq j\leq 6,j\in N)$,  $z_{j}$ with $(j=7,8)$ and $h$.
In this case, we divide the frequency  into
$$|\xi_{l}|\geq {\rm max}\left\{|\xi_{j}|,1\leq j\leq 6,j\neq l,j\in N\right\}(1\leq l\leq 6,l\in N).$$
Without loss of generality, we can assume that
$|\xi_{1}|\geq |\xi_{2}|\geq\cdot\cdot\cdot\geq |\xi_{6}|$ and $|\xi_{7}|\geq |\xi_{8}|$.

\noindent We define $$I_{4}=\left|\int_{\SR}\int_{\SR}J^{\sigma}\partial_{x}
\left(\left(\prod_{j=1}^{6}v_{j}\right)z_{7}z_{8}\right)h(x,t)dxdt\right|.$$
\noindent When $|\xi_{1}|\geq80 |\xi_{2}|,$ $|\xi_{1}|\geq \frac{|\xi_{7}|}{80},$
by using the H\"older inequality,  (\ref{2.04})-(\ref{2.05}), (\ref{2.016}), Lemma 2.3,
we have that
\begin{eqnarray*}
&&I_{4}\leq CN_{1}^{-\frac{1}{10}}
\left\|I^{1/2}I^{1/2}_{-}(J ^{\sigma}v_{1},J ^{\sigma}v_{2})\right\|_{L_{xt}^{2}}
\left(\prod_{j=4}^{6}\|J^{-\frac{2}{7}}v_{j}\|_{L_{xt}^{\infty}}\right)\nonumber\\&&\qquad\times
\left(\prod_{j=7}^{8}\|z_{j}\|_{L_{xt}^{\infty}}\right)
\|I^{\frac{1-\epsilon}{2}}I^{\frac{1-\epsilon}{2}}_{-}(J ^{\sigma}v_{3},h)\|_{L_{xt}^{2}}\leq
CR^{2}\left(\prod_{j=1}^{6}\|v_{j}\|_{X_{\sigma,b}}\right)\|h\|_{X_{0,\frac{1}{2}-\frac{\epsilon}{12}}},
\end{eqnarray*}
outside a set of probability at most
$ C{\rm exp}\left(-C^{'}\frac{R^{2}}{\|\phi\|_{H^{\epsilon}}^{2}}\right).
$

\noindent When $\frac{|\xi_{7}|}{80}\geq|\xi_{1}|\geq80 |\xi_{2}|,$  $|\xi_{7}|\geq 80|\xi_{8}|,$
by using the H\"older inequality,  (\ref{2.04})-(\ref{2.05}), (\ref{2.016}), Lemmas 2.3, 2.4, 2.1,
we have that
\begin{eqnarray*}
&&I_{4}\leq CN_{7}^{-\frac{1}{10}}
\left\|I^{1/2}I^{1/2}_{-}(J ^{s}z_{7},J ^{\sigma}v_{1})\right\|_{L_{xt}^{2}}
\left(\prod_{j=3}^{6}\|J^{-\frac{2}{7}}v_{j}\|_{L_{xt}^{\infty}}\right)
\|z_{8}\|_{L_{xt}^{\infty}}\|I^{\frac{1-\epsilon}{2}}
I^{\frac{1-\epsilon}{2}}_{-}(J ^{\sigma}v_{2},h)\|_{L_{xt}^{2}}\nonumber\\&&\leq
CN_{1}^{-\frac{1}{10}}\left(\prod_{j=1}^{6}\|v_{j}\|_{X_{\sigma,b}}\right)
\|z_{7}\|_{X_{s,c}}\|z_{8}\|_{L_{xt}^{\infty}}
\|h\|_{X_{0,\frac{1}{2}-\frac{\epsilon}{12}}}\nonumber\\
&&\leq CN_{1}^{-\frac{1}{10}}T^{-\frac{\epsilon}{100}}\left(\prod_{j=1}^{6}\|v_{j}\|_{X_{\sigma,b}}\right)
\|P_{N_{7}}\phi^{\omega}\|_{H^{s}}\|z_{8}\|_{L_{xt}^{\infty}}\|h\|_{X_{0,\frac{1}{2}-\frac{\epsilon}{12}}}\nonumber\\
&&\leq CT^{-\frac{\epsilon}{100}}R^{2}
\left(\prod_{j=1}^{6}\|v_{j}\|_{X_{\sigma,b}}\right)\|h\|_{X_{0,\frac{1}{2}-\frac{\epsilon}{12}}},
\end{eqnarray*}
outside a set of probability at most
$
C{\rm exp}\left(-C^{'}\frac{R^{2}}{\|\phi\|_{H^{s}}^{2}}\right).
$

\noindent When $\frac{|\xi_{7}|}{80}\geq|\xi_{1}|\geq80 |\xi_{2}|,$  $|\xi_{7}|\leq 80|\xi_{8}|,$
by using the H\"older inequality,  (\ref{2.05})-(\ref{2.07}),  (\ref{2.016}),  Lemmas 2.4, 2.1, 2.2,
we have that
\begin{eqnarray*}
&&I_{4}\leq CN_{7}^{-\frac{1}{10}}
\left\|I^{1/2}I^{1/2}_{-}(J ^{s}z_{7},J ^{\sigma}v_{1})\right\|_{L_{xt}^{2}}
\left(\prod_{j=2}^{6}\|v_{j}\|_{L_{xt}^{\frac{28}{2-7\epsilon}}}\right)
\|J^{s}z_{8}\|_{L_{xt}^{\frac{56}{1+63\epsilon}}}\|h\|_{L_{xt}^{\frac{8}{1+\epsilon}}}\nonumber\\&&\leq
CN_{1}^{-\frac{1}{10}}\left(\prod_{j=1}^{6}\|v_{j}\|_{X_{\sigma,b}}\right)
\|z_{7}\|_{X_{s,c}}\|J^{s}z_{8}\|_{L_{xt}^{\frac{56}{1+63\epsilon}}}
\|h\|_{X_{0,\frac{1}{2}-\frac{\epsilon}{12}}}\nonumber\\
&&\leq CN_{1}^{-\frac{1}{10}}T^{-\frac{\epsilon}{100}}\left(\prod_{j=1}^{6}\|v_{j}\|_{X_{\sigma,b}}\right)
\|P_{N_{7}}\phi^{\omega}\|_{H^{s}}\|J^{s}z_{8}\|_{L_{xt}^{\frac{56}{1+63\epsilon}}}
\|h\|_{X_{0,\frac{1}{2}-\frac{\epsilon}{12}}}\nonumber\\
&&\leq CT^{-\frac{\epsilon}{100}}R^{2}\left(\prod_{j=1}^{6}\|v_{j}\|_{X_{\sigma,b}}\right)
\|h\|_{X_{0,\frac{1}{2}-\frac{\epsilon}{12}}},
\end{eqnarray*}
outside a set of probability at most
$
\leq C{\rm exp}\left(-C^{'}\frac{R^{2}}{\|\phi\|_{H^{s}}^{2}}\right).
$

\noindent When $|\xi_{1}|\sim |\xi_{2}|\geq 80|\xi_{3}|,$ $|\xi_{1}|\geq \frac{|\xi_{7}|}{80},$
by using the H\"older inequality,   (\ref{2.05})-(\ref{2.07}), (\ref{2.014}),  Lemma 2.2,
we obtain that
\begin{eqnarray*}
&&I_{4}\leq CN_{1}^{-\frac{1}{20}}\|J^{\sigma}v_{1}\|_{L_{xt}^{8}}
\left\|I^{1/2}I^{1/2}_{-}(J^{\sigma}v_{2},J^{\sigma}v_{3})\right\|_{L_{xt}^{2}}
\left(\prod_{j=4}^{6}\|v_{j}\|_{L_{xt}^{\frac{28}{2-7\epsilon}}}\right)
\left(\prod_{j=7}^{8}\|z_{j}\|_{L_{xt}^{\frac{112}{2+35\epsilon}}}\right)\nonumber\\&&\qquad\times
\|h\|_{L_{xt}^{\frac{8}{1+\epsilon}}}\leq
CN_{1}^{-\frac{1}{20}}\left(\prod_{j=1}^{6}\|v_{j}\|_{X_{\sigma,b}}\right)
\left(\prod_{j=7}^{8}\|z_{j}\|_{L_{xt}^{\frac{112}{2+35\epsilon}}}\right)
\|h\|_{X_{0,\frac{1}{2}-\frac{\epsilon}{12}}}\nonumber\\&&
\leq CR^{2}\left(\prod_{j=1}^{6}\|v_{j}\|_{X_{\sigma,b}}\right)
\|h\|_{X_{0,\frac{1}{2}-\frac{\epsilon}{12}}},
\end{eqnarray*}
outside a set of probability at most
$
C{\rm exp}\left(-C^{'}\frac{R^{2}}{T^{\frac{2+35\epsilon}{56}}\|\phi\|_{H^{\epsilon}}^{2}}\right).
$

\noindent When $\frac{|\xi_{7}|}{80}\geq|\xi_{1}|\sim |\xi_{2}|\geq 80|\xi_{3}|,$ $|\xi_{7}|\geq 80|\xi_{8}|,$
this case can be proved similarly to case $\frac{|\xi_{7}|}{80}\geq|\xi_{1}|\geq 80|\xi_{2}|,$ $|\xi_{7}|\geq 80|\xi_{8}|$ of Lemma 3.4.

\noindent When $\frac{|\xi_{7}|}{80}\geq|\xi_{1}|\sim |\xi_{2}|\geq 80|\xi_{3}|,$ $|\xi_{7}|\leq 80|\xi_{8}|,$
this case can be proved similarly to case $\frac{|\xi_{7}|}{80}\geq|\xi_{1}|\geq 80|\xi_{2}|,$ $|\xi_{7}|\leq 80|\xi_{8}|$ of Lemma 3.4.

\noindent When $|\xi_{1}|\sim |\xi_{l}|\geq 80|\xi_{l+1}|(l=3,4),$
$|\xi_{1}|\geq \frac{|\xi_{7}|}{80},$   this case can be proved similarly to case
$|\xi_{1}|\sim |\xi_{2}|\geq 80|\xi_{3}|,$
$|\xi_{1}|\geq \frac{|\xi_{7}|}{80}$ of Lemma 3.4.

\noindent When $\frac{|\xi_{7}|}{80}\geq|\xi_{1}|\sim |\xi_{l}|\geq 80|\xi_{l+1}|(l=3,4),$ $|\xi_{7}|\geq 80|\xi_{8}|,$
this case can be proved similarly to case $\frac{|\xi_{7}|}{80}\geq|\xi_{1}|\geq 80|\xi_{2}|,$ $|\xi_{7}|\geq 80|\xi_{8}|$ of Lemma 3.4.

\noindent When $\frac{|\xi_{7}|}{80}\geq|\xi_{1}|\sim |\xi_{l}|\geq 80|\xi_{l+1}|(l=3,4),$ $|\xi_{7}|\leq 80|\xi_{8}|,$
this case can be proved similarly to case $\frac{|\xi_{7}|}{80}\geq|\xi_{1}|\geq 80|\xi_{2}|,$ $|\xi_{7}|\leq 80|\xi_{8}|$ of Lemma 3.4.

\noindent When $|\xi_{1}|\sim |\xi_{5}|\geq80|\xi_{6}|,$ $|\xi_{1}|\geq\frac{|\xi_{7}|}{80},$
 by using the H\"older inequality, (\ref{2.05})-(\ref{2.07}) and Lemma 2.2,  we have that
\begin{eqnarray*}
&&I_{4}\leq C
\left(\prod_{j=1}^{4}\|v_{j}\|_{L_{xt}^{\frac{28}{2-7\epsilon}}}\right)
\left\|I^{1/2}I^{1/2}_{-}(J^{\sigma}v_{5},J^{\sigma}v_{6})\right\|_{L_{xt}^{2}}
\left(\prod_{j=7}^{8}\|z_{j}\|_{L_{xt}^{\frac{112}{5+49\epsilon}}}\right)
\|h\|_{L_{xt}^{\frac{8}{1+\epsilon}}}\nonumber\\&&\leq CN_{1}^{-2\epsilon}
\left(\prod_{j=1}^{6}\|v_{j}\|_{X_{\sigma,b}}\right)
\left(\prod_{j=7}^{8}\|z_{j}\|_{L_{xt}^{\frac{112}{5+49\epsilon}}}\right)
\|h\|_{X_{0,\frac{1}{2}-\frac{\epsilon}{12}}}\nonumber\\
&&\leq CR^{2}\left(\prod_{j=1}^{6}\|v_{j}\|_{X_{\sigma,b}}\right)
\|h\|_{X_{0,\frac{1}{2}-\frac{\epsilon}{12}}},
\end{eqnarray*}
outside a set of probability at most
$
C{\rm exp}\left(-C^{'}\frac{R^{2}}{T^{\frac{5+49\epsilon }{56}}\|\phi\|_{H^{\epsilon}}^{2}}\right).
$

\noindent When $\frac{|\xi_{7}|}{80}\geq|\xi_{1}|\sim |\xi_{5}|\geq80|\xi_{6}|,$ $|\xi_{7}|\geq 80|\xi_{8}|,$
this case can be proved similarly to case $\frac{|\xi_{7}|}{80}\geq|\xi_{1}|\geq 80|\xi_{2}|,$ $|\xi_{7}|\geq 80|\xi_{8}|$ of Lemma 3.4.

\noindent When $\frac{|\xi_{7}|}{80}\geq|\xi_{1}|\sim |\xi_{5}|\geq80|\xi_{6}|,$ $|\xi_{7}|\leq 80|\xi_{8}|,$
this case can be proved similarly to case $\frac{|\xi_{7}|}{80}\geq|\xi_{1}|\geq 80|\xi_{2}|,$ $|\xi_{7}|\leq 80|\xi_{8}|$ of Lemma 3.4.

\noindent When $|\xi_{1}|\sim |\xi_{6}|,$ $|\xi_{1}|\geq\frac{|\xi_{7}|}{80},$
by using the H\"older inequality, (\ref{2.06}),
(\ref{2.014}) and   Lemma 2.2,  we have that
\begin{eqnarray*}
I_{4}\leq C
\left(\prod_{j=1}^{6}\|v_{j}\|_{X_{\sigma,b}}\right)
\left(\prod_{j=7}^{8}\|z_{j}\|_{L_{xt}^{\frac{16}{1-\epsilon}}}\right)
\|h\|_{X_{0,\frac{1}{2}-\frac{\epsilon}{12}}}\leq CR^{2}
\left(\prod_{j=1}^{6}\|v_{j}\|_{X_{\sigma,b}}\right)\|h\|_{X_{0,\frac{1}{2}-\frac{\epsilon}{12}}},
\end{eqnarray*}
outside a set of probability at most
$
 C{\rm exp}\left(-C^{'}\frac{R^{2}}{T^{\frac{1-\epsilon }{8}}\|\phi\|_{H^{\epsilon}}^{2}}\right).
$

\noindent When $\frac{|\xi_{7}|}{80}\geq|\xi_{1}|\sim |\xi_{6}|,$ $|\xi_{7}|\geq 80|\xi_{8}|,$
this case can be proved similarly to case $\frac{|\xi_{7}|}{80}\geq|\xi_{1}|\geq 80|\xi_{2}|,$
$|\xi_{7}|\geq 80|\xi_{8}|$ of Lemma 3.4.

\noindent When $\frac{|\xi_{7}|}{80}\geq|\xi_{1}|\sim |\xi_{6}|,$ $|\xi_{7}|\leq 80|\xi_{8}|,$
this case can be proved similarly to case $\frac{|\xi_{7}|}{80}\geq|\xi_{1}|\geq 80|\xi_{2}|,$
$|\xi_{7}|\leq 80|\xi_{8}|$ of Lemma 3.4.

\noindent When $|\xi_{1}|\sim |\xi_{7}|,$ $|\xi_{7}|\geq80|\xi_{8}|,$ this case can be proved similarly to case
$|\xi_{1}|\sim |\xi_{6}|,$ $|\xi_{1}|\geq\frac{|\xi_{7}|}{80}$ of Lemma 3.4.

\noindent When $|\xi_{1}|\sim |\xi_{7}|,$ $|\xi_{7}|\leq80|\xi_{8}|,$ this case can be proved similarly to case
$|\xi_{1}|\sim |\xi_{6}|,$ $|\xi_{1}|\geq\frac{|\xi_{7}|}{80}$ of Lemma 3.4.

This completes the proof of Lemma 3.4.

\begin{Lemma}\label{Lemma3.5}
Let $s\geq\frac{17}{112}+\epsilon$ and $\sigma=\frac{3}{14}+2\epsilon$,
$v_{j}=\eta v$ with $(1\leq j\leq 5, j\in N)$  and $z_{j}=\eta_{T}S(t)
\phi^{\omega}(6\leq j\leq8,j\in N)$. Then, we have that
\begin{eqnarray*}
\left|\int_{\SR}\int_{\SR}J^{\sigma}\partial_{x}\left(\left(\prod_{j=1}^{5}v_{j}\right)
\left(\prod_{j=6}^{8}z_{j}\right)\right)h(x,t)dxdt\right|\leq CT^{-\frac{\epsilon}{100}}R^{3}\left(\prod_{j=1}^{5}\|v_{j}\|_{X_{\sigma,b}}\right)\|h\|_{X_{0,\frac{1}{2}-\frac{\epsilon}{12}}},
\end{eqnarray*}
outside a set of probability at most
$
C{\rm exp}\left(-C^{'}\frac{R^{2}}{\|\phi\|_{H^{s}}^{2}}\right).
$

\end{Lemma}
\noindent {\bf Proof.}
We dyadically decompose $v_{j}=\eta v$ with $(1\leq j\leq 5, j\in N)$ and
$z_{j}=\eta_{T}S(t)\phi^{\omega}$
  with $(6\leq j\leq8,j\in Z)$  and $h$ such that frequency supports are
$\left\{|\xi_{j}|\sim N_{j}\right\}$ for some dyadically $N_{j}\geq 1$ and
we still denote them by  $v_{j}=\eta v$ with $(1\leq j\leq 5, j\in N)$,
$z_{j}$ with $(6\leq j\leq 8,j\in N)$ and $h$.
In this case, we divide the frequency  into
$$|\xi_{l}|\geq {\rm max}\left\{|\xi_{j}|,1\leq j\leq 5,j\neq l,j\in N\right\}(1\leq l\leq 5,l\in N).$$
Without loss of generality, we can assume that $|\xi_{1}|\geq |\xi_{2}|
\geq|\xi_{3}|\geq |\xi_{4}|\geq|\xi_{5}|$ and $|\xi_{6}|\geq |\xi_{7}|\geq |\xi_{8}|.$

\noindent We define $$I_{5}=\left|\int_{\SR}\int_{\SR}J^{\sigma}\partial_{x}
\left(\left(\prod_{j=1}^{5}v_{j}\right)
\left(\prod_{j=6}^{8}z_{j}\right)\right)h(x,t)dxdt\right|.$$
\noindent When $|\xi_{1}|\geq80 |\xi_{2}|,$ $|\xi_{1}|\geq 80|\xi_{6}|,$
by using the H\"older inequality, (\ref{2.04})-(\ref{2.05}),
(\ref{2.016}),  Lemma 2.3, we have that
\begin{eqnarray*}
&&I_{5}\leq CN_{1}^{-\frac{1}{2}}
\left\|I^{1/2}I^{1/2}_{-}(J ^{\sigma}v_{1},J ^{\sigma}v_{2})\right\|_{L_{xt}^{2}}
\left(\prod_{j=4}^{5}\|J^{-\frac{2}{7}}v_{j}\|_{L_{xt}^{\infty}}\right)
\left(\prod_{j=6}^{8}\|z_{j}\|_{L_{xt}^{\infty}}\right)\nonumber\\&&\qquad \times
\|I^{\frac{1-\epsilon}{2}}I^{\frac{1-\epsilon}{2}}_{-}(J ^{\sigma}v_{3},h)\|_{L_{xt}^{2}}\nonumber\\&&\leq
CN_{1}^{-\frac{1}{2}}\left(\prod_{j=1}^{5}\|v_{j}\|_{X_{\sigma,b}}\right)
\left(\prod_{j=6}^{8}\|z_{j}\|_{L_{xt}^{\infty}}\right)
\|h\|_{X_{0,\frac{1}{2}-\frac{\epsilon}{12}}}\leq CR^{3}\left(\prod_{j=1}^{5}\|v_{j}\|_{X_{\sigma,b}}\right)
\|h\|_{X_{0,\frac{1}{2}-\frac{\epsilon}{12}}},
\end{eqnarray*}
outside a set of probability at most
$
C{\rm exp}\left(-C^{'}\frac{R^{2}}{\|\phi\|_{H^{\epsilon}}^{2}}\right).
$

\noindent When $|\xi_{1}|\geq80 |\xi_{2}|,$ $\frac{|\xi_{6}|}{80}\leq|\xi_{1}|\leq 80|\xi_{6}|,$
by using the H\"older inequality, (\ref{2.05})-(\ref{2.08}), Lemma 2.2,
we have that
\begin{eqnarray*}
&&I_{5}
\leq CN_{1}^{-\frac{1}{20}}\left\|I^{1/2}I^{1/2}_{-}(J^{\sigma}v_{1},J^{\sigma}v_{2})\right\|_{L_{xt}^{2}}
\|v_{3}\|_{L_{xt}^{\frac{28}{2-7\epsilon}}}\nonumber\\&&\qquad\times
\left(\prod_{j=4}^{5}\|J^{-\frac{\sigma}{5}}v_{j}\|_{L_{xt}^{\frac{280}{17+7\epsilon}}}\right)
\left(\prod_{j=6}^{8}\|J^{s}z_{j}\|_{L_{xt}^{\frac{280}{17+7\epsilon}}}\right)
\|h\|_{L_{xt}^{\frac{8}{1+\epsilon}}}\nonumber\\&&\leq
CN_{1}^{-\frac{1}{20}}\left(\prod_{j=1}^{5}\|v_{j}\|_{X_{\sigma,b}}\right)\left(\prod_{j=6}^{8}
\|z_{j}\|_{L_{xt}^{\frac{280}{17+7\epsilon}}}\right)\|h\|_{X_{0,\frac{1}{2}-\frac{\epsilon}{12}}}\nonumber\\
&&\leq CN_{1}^{-\frac{1}{20}}R^{3}\left(\prod_{j=1}^{5}\|v_{j}\|_{X_{\sigma,b}}\right)
\|h\|_{X_{0,\frac{1}{2}-\frac{\epsilon}{12}}},
\end{eqnarray*}
outside a set of probability at most
$
C{\rm exp}\left(-C^{'}\frac{R^{2}}{T^{\frac{17+7\epsilon}{140}}\|\phi\|_{H^{\epsilon}}^{2}}\right).
$

\noindent When $\frac{|\xi_{6}|}{80}\geq|\xi_{1}|\geq80 |\xi_{2}|,$ $|\xi_{6}|\geq 80|\xi_{7}|,$
by using the H\"older inequality,
(\ref{2.04})-(\ref{2.05}) and (\ref{2.016}), Lemmas 2.4, 2.3,  2.1,
we have that
\begin{eqnarray*}
&&I_{5}\leq CN_{6}^{-\frac{1}{10}}
\left\|I^{1/2}I^{1/2}_{-}(J ^{s}z_{6},J ^{\sigma}v_{1})\right\|_{L_{xt}^{2}}
\left(\prod_{j=3}^{5}\|J^{-\frac{2}{7}}v_{j}\|_{L_{xt}^{\infty}}\right)
\left(\prod_{j=7}^{8}\|z_{j}\|_{L_{xt}^{\infty}}\right)\nonumber\\&&\qquad\times
\|I^{\frac{1-\epsilon}{2}}I^{\frac{1-\epsilon}{2}}_{-}(J ^{\sigma}v_{2},h)\|_{L_{xt}^{2}}\nonumber\\&&\leq
CN_{6}^{-\frac{1}{10}}\left(\prod_{j=1}^{5}\|v_{j}\|_{X_{\sigma,b}}\right)
\|z_{6}\|_{X_{s,c}}
\left(\prod_{j=7}^{8}\|z_{j}\|_{L_{xt}^{\infty}}\right)\|h\|_{X_{0,\frac{1}{2}-\frac{\epsilon}{12}}}\nonumber\\
&&\leq CN_{6}^{-\frac{1}{10}}T^{-\frac{\epsilon}{100}}
\left(\prod_{j=1}^{5}\|v_{j}\|_{X_{\sigma,b}}\right)\|P_{N_{6}}\phi^{\omega}\|_{H^{s}}
\left(\prod_{j=7}^{8}\|z_{j}\|_{L_{xt}^{\infty}}\right)
\|h\|_{X_{0,\frac{1}{2}-\frac{\epsilon}{12}}}\nonumber\\&&\leq CT^{-\frac{\epsilon}{100}}
R^{3}\left(\prod_{j=1}^{5}\|v_{j}\|_{X_{\sigma,b}}\right)
\|h\|_{X_{0,\frac{1}{2}-\frac{\epsilon}{12}}},
\end{eqnarray*}
outside a set of probability at most
$
 C{\rm exp}\left(-C^{'}\frac{R^{2}}{\|\phi\|_{H^{s}}^{2}}\right).
$

\noindent When $\frac{|\xi_{6}|}{80}\geq|\xi_{1}|\geq80 |\xi_{2}|,$
$|\xi_{6}|\leq 80|\xi_{7}|,$  by using the H\"older inequality,
(\ref{2.05})-(\ref{2.06}), (\ref{2.015}),  Lemmas 2.4, 2.1,  2.3,
we have that
\begin{eqnarray*}
&&I_{5}\leq CN_{6}^{-\frac{1}{40}}
\left(\prod_{j=2}^{5}\|J^{\sigma-\frac{1+\epsilon}{8}}v_{j}\|_{L_{xt}^{\frac{32}{3-\epsilon}}}\right)
\left\|I^{1/2}I^{1/2}_{-}(J ^{s}z_{6},J ^{\sigma}v_{1})\right\|_{L_{xt}^{2}}
\left(\prod_{j=7}^{8}\left\|J^{s}z_{j}\right\|_{L_{xt}^{\infty}}\right)
\|h\|_{L_{xt}^{\frac{8}{1+\epsilon}}}
\nonumber\\&&\leq
CN_{6}^{-\frac{1}{40}}\left(\prod_{j=1}^{5}\|v_{j}\|_{X_{\sigma,b}}\right)
\|z_{6}\|_{X_{s,c}}
\left(\prod_{j=7}^{8}\|J^{s}z_{j}\|_{L_{xt}^{\infty}}\right)
\|h\|_{X_{0,\frac{1}{2}-\frac{\epsilon}{12}}}\nonumber\\
&&\leq CN_{6}^{-\frac{1}{40}}T^{-\frac{\epsilon}{100}}
\left(\prod_{j=1}^{5}\|v_{j}\|_{X_{\sigma,b}}\right)
\|P_{N_{6}}\phi^{\omega}\|_{H^{s}}
\left(\prod_{j=7}^{8}\|J^{s}z_{j}\|_{L_{xt}^{\infty}}\right)
\|h\|_{X_{0,\frac{1}{2}-\frac{\epsilon}{12}}}\nonumber\\&&\leq CT^{-\frac{\epsilon}{100}}
R^{3}\left(\prod_{j=1}^{5}\|v_{j}\|_{X_{\sigma,b}}\right)
\|h\|_{X_{0,\frac{1}{2}-\frac{\epsilon}{12}}},
\end{eqnarray*}
outside a set of probability at most
$
C{\rm exp}\left(-C^{'}\frac{R^{2}}{\|\phi\|_{H^{s}}^{2}}\right).
$

\noindent When $|\xi_{1}|\sim |\xi_{2}|\geq 80|\xi_{3}|,$ $|\xi_{1}|\geq \frac{|\xi_{6}|}{80},$
 by using the H\"older inequality,     (\ref{2.05})-(\ref{2.08}),
 Lemma 2.2,
we have that
\begin{eqnarray*}
&&I_{5}\leq CN_{1}^{-\frac{1}{20}}\|J^{s}v_{1}\|_{L_{xt}^{8}}
\left\|I^{1/2}I^{1/2}_{-}(J^{\sigma}v_{2},J^{\sigma}v_{3})\right\|_{L_{xt}^{2}}
\left(\prod_{j=4}^{5}\|v_{j}\|_{L_{xt}^{\frac{28}{2-7\epsilon}}}\right)
\left(\prod_{j=6}^{8}\|z_{j}\|_{L_{xt}^{\frac{168}{6+21\epsilon}}}\right)\nonumber\\&&\qquad \times
\|h\|_{L_{xt}^{\frac{8}{1+\epsilon}}}\nonumber\\&&\leq
CN_{1}^{-\frac{1}{20}}N_{1}^{-2\epsilon}
\left(\prod_{j=1}^{5}\|v_{j}\|_{X_{\sigma,b}}\right)\left(\prod_{j=6}^{8}
\|z_{j}\|_{L_{xt}^{\frac{168}{6+21\epsilon}}}\right)
\|h\|_{X_{0,\frac{1}{2}-\frac{\epsilon}{12}}}\nonumber\\
&&\leq CN_{1}^{-\frac{1}{20}}R^{3}
\left(\prod_{j=1}^{5}\|v_{j}\|_{X_{\sigma,b}}\right)\|h\|_{X_{0,\frac{1}{2}-\frac{\epsilon}{12}}},
\end{eqnarray*}
outside a set of probability at most
$
C{\rm exp}\left(-C^{'}\frac{R^{2}}{T^{\frac{6+21\epsilon}{84}}\|\phi\|_{H^{s}}^{2}}\right).
$

\noindent When $ \frac{|\xi_{6}|}{80}\geq|\xi_{1}|\sim |\xi_{2}|\geq 80|\xi_{3}|,$
$|\xi_{6}|\geq 80|\xi_{7}|,$ this case can be proved similarly to case
$ \frac{|\xi_{6}|}{80}\geq|\xi_{1}|\geq 80|\xi_{2}|,$
 $|\xi_{6}|\geq 80|\xi_{7}|$ of Lemma 3.5.

\noindent When $ \frac{|\xi_{6}|}{80}\geq|\xi_{1}|\sim |\xi_{2}|\geq 80|\xi_{3}|,$
 $|\xi_{6}|\leq 80|\xi_{7}|,$ this case can be proved similarly to  case
$\frac{|\xi_{6}|}{80}\geq|\xi_{1}|\geq80 |\xi_{2}|$ of Lemma 3.5.

\noindent When $|\xi_{1}|\sim |\xi_{3}|\geq 80|\xi_{4}|,$
$|\xi_{1}|\geq \frac{|\xi_{6}|}{80},$
this case can be proved similarly to case $|\xi_{1}|\sim |\xi_{2}|\geq 80|\xi_{3}|,$
$|\xi_{1}|\geq \frac{|\xi_{6}|}{80}$ of Lemma 3.5.

\noindent When $\frac{|\xi_{6}|}{80}\geq|\xi_{1}|\sim |\xi_{3}|\geq 80|\xi_{4}|,$
this case can be proved similarly to  case $\frac{|\xi_{6}|}{80}\geq|\xi_{1}|\sim |\xi_{2}|\geq 80|\xi_{3}|$ of Lemma 3.5.

\noindent When $|\xi_{1}|\sim |\xi_{4}|\geq 80|\xi_{5}|,$ $|\xi_{1}|\geq \frac{|\xi_{6}|}{80},$
this case can be proved similarly to case $|\xi_{1}|\sim |\xi_{2}|\geq 80|\xi_{3}|,$
$|\xi_{1}|\geq \frac{|\xi_{6}|}{80}$ of Lemma 3.5.

\noindent When $\frac{|\xi_{6}|}{80}\geq|\xi_{1}|\sim |\xi_{4}|\geq 80|\xi_{5}|,$
 this case can be proved similarly to case $\frac{|\xi_{6}|}{80}\geq|\xi_{1}|\sim |\xi_{2}|\geq 80|\xi_{3}|$ of Lemma 3.5.

\noindent When $|\xi_{1}|\sim |\xi_{5}|, |\xi_{1}|\geq80|\xi_{6}|,$
by using  the H\"older inequality,
(\ref{2.05})-(\ref{2.07}),  Lemmas 2.4, 2.1, 2.2, we have that
\begin{eqnarray*}
&&I_{5}\leq C
\left(\prod_{j=2}^{5}\|v_{j}\|_{L_{xt}^{\frac{28}{2-7\epsilon}}}\right)
\left\|I^{1/2}I^{1/2}_{-}(J^{\sigma}v_{1},J^{s}z_{6})\right\|_{L_{xt}^{2}}
\left(\prod_{j=7}^{8}\|z_{j}\|_{L_{xt}^{\frac{112}{5+49\epsilon}}}\right)
\|h\|_{L_{xt}^{\frac{8}{1+\epsilon}}}\nonumber\\&&\leq CN_{1}^{-2\epsilon}
\left(\prod_{j=1}^{5}\|v_{j}\|_{X_{\sigma,b}}\right)\|z_{6}\|_{X_{s,c}}
\left(\prod_{j=7}^{8}\|z_{j}\|_{L_{xt}^{\frac{112}{5+49\epsilon}}}\right)
\|h\|_{X_{0,\frac{1}{2}-\frac{\epsilon}{12}}}
\nonumber\\
&&\leq CT^{-\frac{\epsilon}{100}}N_{1}^{-2\epsilon}
\left(\prod_{j=1}^{5}\|v_{j}\|_{X_{\sigma,b}}\right)\|P_{N_{6}}\phi^{\omega}\|_{H^{s}}
\left(\prod_{j=7}^{8}\|z_{j}\|_{L_{xt}^{\frac{112}{5+49\epsilon}}}\right)
\|h\|_{X_{0,\frac{1}{2}-\frac{\epsilon}{12}}}\nonumber\\
&&\leq CT^{-\frac{\epsilon}{100}}R^{3}
\left(\prod_{j=1}^{5}\|v_{j}\|_{X_{\sigma,b}}\right)\|h\|_{X_{0,\frac{1}{2}-\frac{\epsilon}{12}}},
\end{eqnarray*}
outside a set of probability at most
$
C{\rm exp}\left(-C^{'}\frac{R^{2}}{\|\phi\|_{H^{s}}^{2}}\right).
$

\noindent When $|\xi_{1}|\sim |\xi_{5}|, \frac{|\xi_{6}|}{80}\leq|\xi_{1}|\leq80|\xi_{6}|,$
 by using the H\"older inequality, (\ref{2.06}),
  (\ref{2.014}), (\ref{2.016}), Lemmas 2.2,
  we have that
\begin{eqnarray*}
&&I_{5}\leq CN_{6}^{-\frac{1}{112}}
\left(\prod_{j=1}^{5}\|J^{\sigma}v_{j}\|_{L_{xt}^{8}}\right)
\left(\prod_{j=6}^{8}\|J^{s}z_{j}\|_{L_{xt}^{\frac{24}{2-\epsilon}}}\right)
\|h\|_{L_{xt}^{\frac{8}{1+\epsilon}}}\nonumber\\
&&\leq CN_{6}^{-\frac{1}{112}}
\left(\prod_{j=1}^{5}\|v_{j}\|_{X_{\sigma,b}}\right)
\left(\prod_{j=6}^{8}\|J^{s}z_{j}\|_{L_{xt}^{\frac{24}{2-\epsilon}}}\right)
\|h\|_{X_{0,\frac{1}{2}-\frac{\epsilon}{12}}}\nonumber\\
&&\leq CN_{6}^{-\frac{1}{112}}
\left(\prod_{j=1}^{5}\|v_{j}\|_{X_{\sigma,b}}\right)
\left(\prod_{j=6}^{8}\|J^{s}z_{j}\|_{L_{xt}^{\frac{24}{2-\epsilon}}}\right)
\|h\|_{X_{0,\frac{1}{2}-\frac{\epsilon}{12}}}\nonumber\\
&&\leq CR^{3}
\left(\prod_{j=1}^{5}\|v_{j}\|_{X_{\sigma,b}}\right)\|h\|_{X_{0,\frac{1}{2}-\frac{\epsilon}{12}}},
\end{eqnarray*}
outside a set of probability at most
$
 C{\rm exp}\left(-C^{'}\frac{R^{2}}{T^{\frac{2-\epsilon}{12}}\|\phi\|_{H^{s}}^{2}}\right).
$

\noindent When $\frac{|\xi_{6}|}{80}\geq|\xi_{1}|\sim |\xi_{5}|, |\xi_{6}|\geq 80|\xi_{7}|,$
this case can be proved similarly to case  $ \frac{|\xi_{6}|}{80}\geq|\xi_{1}|\geq 80|\xi_{2}|,$
$|\xi_{6}|\geq 80|\xi_{7}|$ of Lemma 3.5.

\noindent When $\frac{|\xi_{6}|}{80}\geq|\xi_{1}|\sim |\xi_{5}|, |\xi_{6}|\leq 80|\xi_{7}|,$
this case can be proved similarly to  case
$\frac{|\xi_{6}|}{80}\geq|\xi_{1}|\geq80 |\xi_{2}|$,
$|\xi_{6}|\leq 80|\xi_{7}|$ of Lemma 3.5.

\noindent When $|\xi_{1}|\sim |\xi_{6}|,$  this case can be proved similarly to case
 $|\xi_{1}|\sim |\xi_{5}|, \frac{|\xi_{6}|}{80}\leq|\xi_{1}|\leq80|\xi_{6}|$
of Lemma 3.5.

This completes the proof of Lemma 3.5.

\begin{Lemma}\label{Lemma3.6}
Let $s\geq \frac{17}{112}+\epsilon$ and $\sigma=\frac{3}{14}+2\epsilon$
and $v_{j}=\eta v(1\leq j\leq 4,j\in N)$ and $z_{j}=\eta_{T}S(t)\phi^{\omega}(5\leq j\leq8,j\in N)$.
Then, we have that
\begin{eqnarray*}
\left|\int_{\SR}\int_{\SR}J^{\sigma}\partial_{x}\left(\left(\prod_{j=1}^{4}v_{j}\right)
\left(\prod_{j=5}^{8}z_{j}\right)\right)h(x,t)dxdt\right|\leq C
T^{-\frac{\epsilon}{50}}R^{4}\left(\prod_{j=1}^{4}\|v_{j}\|_{X_{\sigma,b}}\right)
\|h\|_{X_{0,\frac{1}{2}-\frac{\epsilon}{12}}},
\end{eqnarray*}
outside a set of probability at most
$
C{\rm exp}\left(-C^{'}\frac{R^{2}}{\|\phi\|_{H^{s}}^{2}}\right).
$
\end{Lemma}
\noindent {\bf Proof.} We dyadically decompose $v_{j}$ with $(1\leq j\leq 4,j\in N)$
and $z_{j}$ with $(5\leq j\leq8,j\in N)$  and $h$ such that frequency supports are
$\left\{|\xi_{j}|\sim N_{j}\right\}$ for some dyadically $N_{j}\geq 1$
and we still denote them by $v_{j}$ with $(1\leq j\leq 4,j\in N)$ and $z_{j}$ with
$(5\leq j\leq 8,j\in N)$ and $h$.
In this case, we divide the frequency  into
$$|\xi_{l}|\geq {\rm max}\left\{|\xi_{j}|,1\leq j\leq 4,j\neq l,j\in N\right\}(1\leq l\leq 4,l\in N).$$
Without loss of generality, we can assume that
$|\xi_{1}|\geq |\xi_{2}|\geq|\xi_{3}\geq |\xi_{4}| $ and $|\xi_{5}|\geq|\xi_{6}|\geq|\xi_{7}|\geq |\xi_{8}|$.

\noindent We define $$I_{6}=\left|\int_{\SR}\int_{\SR}J^{\sigma}\partial_{x}\left(\left(\prod_{j=1}^{4}v_{j}\right)
\left(\prod_{j=5}^{8}z_{j}\right)\right)h(x,t)dxdt\right|.$$
\noindent When $|\xi_{1}|\geq80 |\xi_{2}|,$ $|\xi_{1}|\geq 80|\xi_{5}|,$
by using the H\"older inequality, (\ref{2.04})-(\ref{2.05}),
 (\ref{2.016})  and   Lemma 2.3, we have that
\begin{eqnarray*}
&&I_{6}\leq CN_{1}^{-\frac{1}{20}}
\left\|I^{1/2}I^{1/2}_{-}(J^{\sigma}v_{1},J ^{\sigma}v_{2})\right\|_{L_{xt}^{2}}
\|I^{\frac{1-\epsilon}{2}}I^{\frac{1-\epsilon}{2}}_{-}(J ^{\sigma}v_{3},h)\|_{L_{xt}^{2}}\|J^{-\frac{2}{7}}v_{4}\|_{L_{xt}^{\infty}}
\left(\prod_{j=5}^{8}\|z_{j}\|_{L_{xt}^{\infty}}\right)\nonumber\\&&
\leq
CN_{1}^{-\frac{1}{20}}\left(\prod_{j=1}^{4}\|v_{j}\|_{X_{\sigma,b}}\right)
\left(\prod_{j=5}^{8}\|z_{j}\|_{L_{xt}^{\infty}}\right)\|h\|_{X_{0,\frac{1}{2}-\frac{\epsilon}{12}}}\leq CR^{4}\left(\prod_{j=1}^{4}\|v_{j}\|_{X_{\sigma,b}}\right)
\|h\|_{X_{0,\frac{1}{2}-\frac{\epsilon}{12}}},
\end{eqnarray*}
outside a set of probability at most
$
C{\rm exp}\left(-C^{'}\frac{R^{2}}{\|\phi\|_{H^{\epsilon}}^{2}}\right).
$

\noindent When $|\xi_{1}|\geq80 |\xi_{2}|,$ $\frac{|\xi_{5}|}{80}\leq|\xi_{1}|\leq 80|\xi_{5}|,$
by using the H\"older inequality, (\ref{2.05})-(\ref{2.06}), (\ref{2.014})  and  Lemma 2.2,
  we have that
\begin{eqnarray*}
&&I_{6}\leq C
N_{1}^{-\frac{1}{16}}\left\|I^{1/2}I^{1/2}_{-}(J ^{\sigma}v_{1},J ^{\sigma}v_{2})\right\|_{L_{xt}^{2}}
\left(\prod_{j=3}^{4}\|v_{j}\|_{L_{xt}^{8}}\right)
\left(\prod_{j=5}^{8}\|J^{s}z_{j}\|_{L_{xt}^{\frac{32}{1-\epsilon}}}\right)
\|h\|_{L_{xt}^{\frac{8}{1+\epsilon}}}\nonumber\\&&\leq
CN_{1}^{-\frac{1}{16}}\left(\prod_{j=1}^{4}\|v_{j}\|_{X_{\sigma,b}}\right)
\left(\prod_{j=5}^{8}\|J^{s}z_{j}\|_{L_{xt}^{\frac{32}{1-\epsilon}}}\right)
\|h\|_{X_{0,\frac{1}{2}-\frac{\epsilon}{12}}}\nonumber\\
&&\leq CR^{4}\left(\prod_{j=1}^{4}\|v_{j}\|_{X_{\sigma,b}}\right)
\|h\|_{X_{0,\frac{1}{2}-\frac{\epsilon}{12}}},
\end{eqnarray*}
outside a set of probability at most
$
 C{\rm exp}\left(-C^{'}\frac{R^{2}}{T^{\frac{1-\epsilon}{16}}\|\phi\|_{H^{s}}^{2}}\right).
$

\noindent When $\frac{|\xi_{5}|}{80}\geq|\xi_{1}|\geq80 |\xi_{2}|,$ $|\xi_{5}|\geq 80|\xi_{6}|,$
by using the H\"older inequality, (\ref{2.04})-(\ref{2.05}), (\ref{2.016})
 and  Lemmas 2.4, 2.1, 2.3,  we have that
\begin{eqnarray*}
&&I_{6}\leq CN_{5}^{-\frac{1}{20}}
\left\|I^{1/2}I^{1/2}_{-}(J ^{\sigma}v_{1},J ^{s}z_{5})\right\|_{L_{xt}^{2}}\left(\prod_{j=3}^{4}
\|J^{-\frac{2}{7}}v_{j}\|_{L_{xt}^{\infty}}\right)\nonumber\\&&\qquad \times
\left(\prod_{j=6}^{8}\|z_{j}\|_{L_{xt}^{\infty}}\right)
\|I^{\frac{1-\epsilon}{2}}I^{\frac{1-\epsilon}{2}}_{-}(J ^{\sigma}v_{2},h)\|_{L_{xt}^{2}}\nonumber\\&&\leq
CN_{5}^{-\frac{1}{20}}\left(\prod_{j=1}^{4}\|v_{j}\|_{X_{\sigma,b}}\right)
\|z_{5}\|_{X_{s,c}}\left(\prod_{j=6}^{8}\|z_{j}\|_{L_{xt}^{\infty}}\right)
\|h\|_{X_{0,\frac{1}{2}-\frac{\epsilon}{12}}}\nonumber\\
&&\leq CN_{5}^{-\frac{1}{20}}T^{-\frac{\epsilon}{100}}
\left(\prod_{j=1}^{4}\|v_{j}\|_{X_{\sigma,b}}\right)\|P_{N_{5}}\phi^{\omega}\|_{H^{s}}
\left(\prod_{j=6}^{8}\|z_{j}\|_{L_{xt}^{\infty}}\right)
\|h\|_{X_{0,\frac{1}{2}-\frac{\epsilon}{12}}}\nonumber\\
&&\leq CT^{-\frac{\epsilon}{100}}R^{4}\left(\prod_{j=1}^{4}\|v_{j}\|_{X_{\sigma,b}}\right)
\|h\|_{X_{0,\frac{1}{2}-\frac{\epsilon}{12}}},
\end{eqnarray*}
outside a set of probability at most
$
 C{\rm exp}\left(-C^{'}\frac{R^{2}}{\|\phi\|_{H^{s}}^{2}}\right).
$

\noindent When $\frac{|\xi_{5}|}{80}\geq|\xi_{1}|\geq80 |\xi_{2}|,$ $|\xi_{5}|\leq 80|\xi_{6}|,$
by using the H\"older inequality, (\ref{2.05})-(\ref{2.07}),
 Lemmas 2.4, 2.1, 2.2,   we infer that
\begin{eqnarray*}
&&I_{6}\leq CN_{5}^{-\frac{1}{20}}
\left\|I^{1/2}I^{1/2}_{-}(J ^{\sigma}v_{1},J ^{s}z_{5})\right\|_{L_{xt}^{2}}
\left(\prod_{j=2}^{4}\|v_{j}\|_{L_{xt}^{\frac{28}{2-7\epsilon}}}\right)
\left(\prod_{j=6}^{8}\|J^{s}z_{j}\|_{L_{xt}^{\frac{168}{9+35\epsilon}}}\right)
\|h\|_{L_{xt}^{\frac{8}{1+\epsilon}}}\nonumber\\&&\leq
CN_{5}^{-\frac{1}{20}}\left(\prod_{j=1}^{4}\|v_{j}\|_{X_{\sigma,b}}\right)
\|z_{5}\|_{X_{s,c}}\left(\prod_{j=6}^{8}\|z_{j}\|_{L_{xt}^{\frac{168}{9+35\epsilon}}}\right)
\|h\|_{X_{0,\frac{1}{2}-\frac{\epsilon}{12}}}\nonumber\\
&&\leq CT^{-\frac{\epsilon}{100}}N_{5}^{-\frac{1}{20}}
\left(\prod_{j=1}^{4}\|v_{j}\|_{X_{\sigma,b}}\right)\|P_{N_{5}}\phi^{\omega}\|_{H^{s}}
\left(\prod_{j=6}^{8}\|z_{j}\|_{L_{xt}^{\frac{168}{9+35\epsilon}}}\right)
\|h\|_{X_{0,\frac{1}{2}-\frac{\epsilon}{12}}}\nonumber\\&&\leq
CT^{-\frac{\epsilon}{100}}R^{4}\left(\prod_{j=1}^{4}\|v_{j}\|_{X_{\sigma ,b}}\right)
\|h\|_{X_{0,\frac{1}{2}-\frac{\epsilon}{12}}},
\end{eqnarray*}
outside a set of probability at most
$
C{\rm exp}\left(-C^{'}\frac{R^{2}}{\|\phi\|_{H^{s}}^{2}}\right).
$

\noindent When $|\xi_{1}|\sim |\xi_{2}|\geq 80|\xi_{3}|,$ $|\xi_{1}|\geq\frac{|\xi_{5}|}{80},$
by using the H\"older inequality, (\ref{2.05})-(\ref{2.07}),
 Lemma 2.3,
we have that
\begin{eqnarray*}
&&I_{6}\leq CN_{1}^{-\sigma}\|J^{\sigma}v_{1}\|_{L_{xt}^{8}}
\left\|I^{1/2}I^{1/2}_{-}(J^{\sigma}v_{2},J^{\sigma}v_{3})\right\|_{L_{xt}^{2}}
\|J^{\sigma}v_{4}\|_{L_{xt}^{8}}\left(\prod_{j=5}^{8}
\|z_{j}\|_{L_{xt}^{\frac{32}{1-\epsilon}}}\right)
\|h\|_{L_{xt}^{\frac{8}{1+\epsilon}}}\nonumber\\&&\leq
CN_{1}^{-\sigma}\left(\prod_{j=1}^{4}\|v_{j}\|_{X_{\sigma,b}}\right)
\left(\prod_{j=5}^{8}\|z_{j}\|_{L_{xt}^{\frac{32}{1-\epsilon}}}\right)
\|h\|_{X_{0,\frac{1}{2}-\frac{\epsilon}{12}}}\nonumber\\
&&\leq CR^{4}\left(\prod_{j=1}^{4}\|v_{j}\|_{X_{\sigma,b}}\right)
\|h\|_{X_{0,\frac{1}{2}-\frac{\epsilon}{12}}},
\end{eqnarray*}
outside a set of probability at most
$
C{\rm exp}\left(-C^{'}\frac{R^{2}}{T^{\frac{1-\epsilon}{16}}\|\phi\|_{H^{\epsilon}}^{2}}\right).
$

\noindent When $\frac{|\xi_{5}|}{80}\geq|\xi_{1}|\sim |\xi_{2}|\geq 80|\xi_{3}|,$ $|\xi_{5}|\geq80|\xi_{6}|,$
this case can be proved similarly to case $\frac{|\xi_{5}|}{80}\geq|\xi_{1}|\geq 80|\xi_{3}|,$ $|\xi_{5}|\geq80|\xi_{6}|$
of  Lemma 3.6.

\noindent When $\frac{|\xi_{5}|}{80}\geq|\xi_{1}|\sim |\xi_{2}|\geq80 |\xi_{3}|,$ $|\xi_{5}|\leq 80|\xi_{6}|,$
this case can be proved similarly to case $\frac{|\xi_{5}|}{80}\geq|\xi_{1}|\geq80 |\xi_{2}|,$ $|\xi_{5}|\leq 80|\xi_{6}|$
of Lemma 3.6.

\noindent When $|\xi_{1}|\sim |\xi_{3}|\geq 80|\xi_{4}|,$ $|\xi_{1}|\geq \frac{|\xi_{5}|}{80},$
 this case can be proved similarly to case $|\xi_{1}|\geq 80|\xi_{2}|,$
 $|\xi_{1}|\geq\frac{|\xi_{5}|}{80}$ of Lemma 3.6.

\noindent When $\frac{|\xi_{5}|}{80}\geq|\xi_{1}|\sim |\xi_{3}|\geq 80|\xi_{4}|,$ $|\xi_{5}|\geq 80|\xi_{6}|,$
 this case can be proved similarly to case $\frac{|\xi_{5}|}{80}\geq|\xi_{1}|\sim |\xi_{2}|\geq 80|\xi_{3}|,$ $|\xi_{5}|\geq 80|\xi_{6}|$
 of  Lemma 3.6.

\noindent When $\frac{|\xi_{5}|}{80}\geq|\xi_{1}|\sim |\xi_{3}|\geq 80|\xi_{4}|,$
$|\xi_{6}|\leq|\xi_{5}|\leq 80|\xi_{6}|,$  this case can be proved similarly to case
 $\frac{|\xi_{5}|}{80}\geq|\xi_{1}|\geq 80|\xi_{2}|,$
$|\xi_{6}|\leq|\xi_{5}|\leq 80|\xi_{6}|$ of Lemma 3.6.

\noindent When $|\xi_{1}|\sim |\xi_{4}|,  |\xi_{1}|\geq 80|\xi_{5}|,$  by using the
H\"older inequality, (\ref{2.05})-(\ref{2.07}),
 Lemma 2.3, we have that
\begin{eqnarray*}
&&I_{6}\leq C
\left(\prod_{j=2}^{4}\|v_{j}\|_{L_{xt}^{\frac{28}{2-7\epsilon}}}\right)
\left\|I^{1/2}I^{1/2}_{-}(J^{\sigma}v_{1},J^{s}z_{5})\right\|_{L_{xt}^{2}}
\left(\prod_{j=6}^{8}\|z_{j}\|_{L_{xt}^{\frac{168}{9+35\epsilon}}}\right)
\|h\|_{L_{xt}^{\frac{8}{1+\epsilon}}}\nonumber\\&&\leq C
\left(\prod_{j=1}^{4}\|v_{j}\|_{X_{\sigma,b}}\right)
\left(\prod_{j=5}^{8}\|z_{j}\|_{L_{xt}^{\frac{168}{9+14\epsilon}}}\right)
\|h\|_{X_{0,\frac{1}{2}-\frac{\epsilon}{12}}}\nonumber\\
&&\leq CR^{4}\left(\prod_{j=1}^{4}\|v_{j}\|_{X_{\sigma,b}}\right)
\|h\|_{X_{0,\frac{1}{2}-\frac{\epsilon}{12}}},
\end{eqnarray*}
outside a set of probability at most
$
C{\rm exp}\left(-C^{'}\frac{R^{2}}{T^{\frac{9+14\epsilon}{84}}\|\phi\|_{H^{s}}^{2}}\right).
$

\noindent When $|\xi_{1}|\sim |\xi_{4}|,\frac{|\xi_{5}|}{80}\leq|\xi_{1}|\leq 80|\xi_{5}|$, $|\xi_{5}|\geq 80|\xi_{6}|,$
by using the H\"older inequality, (\ref{2.04})-(\ref{2.05}),  Lemmas 2.4, 2.1, 2.2,
 and (\ref{2.016}), we have that
\begin{eqnarray*}
&&I_{6}\leq CN_{5}^{-s}
\|J^{\sigma}v_{1}\|_{L_{xt}^{8}}\left(\prod_{j=2}^{4}\|v_{j}\|_{L_{xt}^{\frac{24}{2-\epsilon}}}\right)
\left\|I^{1/2}I^{1/2}_{-}(J ^{s}z_{5},J ^{s}z_{6})\right\|_{L_{xt}^{2}}
\left(\prod_{j=7}^{8}\|z_{j}\|_{L_{xt}^{\infty}}\right)
\|h\|_{L_{xt}^{\frac{8}{1+\epsilon}}}\nonumber\\&&\leq
CN_{5}^{-s}\left(\prod_{j=1}^{4}\|v_{j}\|_{X_{\sigma,b}}\right)
\left(\prod_{j=5}^{6}\|z_{j}\|_{X_{s,c}}\right)
\left(\prod_{j=7}^{8}\|z_{j}\|_{L_{xt}^{\infty}}\right)
\|h\|_{X_{0,\frac{1}{2}-\frac{\epsilon}{12}}}\nonumber\\
&&\leq CN_{5}^{-s}T^{-\frac{\epsilon}{50}}
\left(\prod_{j=1}^{4}\|v_{j}\|_{X_{\sigma,b}}\right)
\left(\prod_{j=5}^{6}\|P_{N_{j}}\phi^{\omega}\|_{H^{s}}\right)
\left(\prod_{j=7}^{8}\|z_{j}\|_{L_{xt}^{\infty}}\right)
\|h\|_{X_{0,\frac{1}{2}-\frac{\epsilon}{12}}}\nonumber\\
&&\leq CT^{-\frac{\epsilon}{50}}R^{4}\left(\prod_{j=1}^{4}\|v_{j}\|_{X_{\sigma,b}}\right)
\|h\|_{X_{0,\frac{1}{2}-\frac{\epsilon}{12}}},
\end{eqnarray*}
outside a set of probability at most
$
C{\rm exp}\left(-C^{'}\frac{R^{2}}{\|\phi\|_{H^{s}}^{2}}\right).
$

\noindent When $|\xi_{1}|\sim |\xi_{4}|,\frac{|\xi_{5}|}{80}\leq|\xi_{1}|\leq 80|\xi_{5}|$,
$|\xi_{5}|\leq 80|\xi_{6}|,|\xi_{6}|\geq80|\xi_{7}|$,  this case can be proved similarly to case
 $|\xi_{1}|\sim |\xi_{4}|,\frac{|\xi_{5}|}{80}\leq|\xi_{1}|\leq 80|\xi_{5}|$,
 $|\xi_{5}|\geq 80|\xi_{6}|$ of Lemma 3.6.

\noindent When $|\xi_{1}|\sim |\xi_{4}|,\frac{|\xi_{5}|}{80}\leq|\xi_{1}|\leq 80|\xi_{5}|$,
$|\xi_{5}|\leq 80|\xi_{6}|,|\xi_{6}|\leq80|\xi_{7}|,|\xi_{7}|\geq 80|\xi_{8}|$,
this case can be proved similarly to case
 $|\xi_{1}|\sim |\xi_{4}|,\frac{|\xi_{5}|}{80}\leq|\xi_{1}|\leq 80|\xi_{5}|$,
 $|\xi_{5}|\geq 80|\xi_{6}|$ of Lemma 3.6.

\noindent
When $|\xi_{1}|\sim |\xi_{4}|,\frac{|\xi_{5}|}{80}\leq|\xi_{1}|\leq 80|\xi_{5}|$,
$|\xi_{5}|\leq 80|\xi_{6}|,|\xi_{6}|\leq80|\xi_{7}|,|\xi_{7}|\leq 80|\xi_{8}|$,
by using the H\"older inequality, (\ref{2.014}), Lemma 2.2,  we have that
\begin{eqnarray*}
&&I_{6}\leq CN_{1}^{-\frac{1}{20}}
\left(\prod_{j=1}^{4}\|J^{\sigma}v_{j}\|_{L_{xt}^{8}}\right)
\left(\prod_{j=5}^{8}\|J^{s}z_{j}\|_{L_{xt}^{\frac{32}{3-\epsilon}}}\right)
\|h\|_{L_{xt}^{\frac{8}{1+\epsilon}}}\nonumber\\&&\leq  CN_{1}^{-\frac{1}{20}}
\left(\prod_{j=1}^{4}\|v_{j}\|_{X_{\sigma,b}}\right)
\left(\prod_{j=5}^{8}\|J^{s}z_{j}\|_{L_{xt}^{\frac{32}{3-\epsilon}}}\right)
\|h\|_{X_{0,\frac{1}{2}-\frac{\epsilon}{12}}}\nonumber\\
&&\leq CR^{4}\left(\prod_{j=1}^{4}\|v_{j}\|_{X_{\sigma,b}}\right)
\|h\|_{X_{0,\frac{1}{2}-\frac{\epsilon}{12}}},
\end{eqnarray*}
outside a set of probability at most
$
C{\rm exp}\left(-C^{'}\frac{R^{2}}{T^{\frac{3-\epsilon}{16}}\|\phi\|_{H^{s}}^{2}}\right).
$

\noindent When $\frac{|\xi_{5}|}{80}\geq|\xi_{1}|\sim |\xi_{4}|,|\xi_{5}|\geq 80|\xi_{6}|$,
this case can be proved similarly to case $\frac{|\xi_{5}|}{80}\geq|\xi_{1}|\geq 80|\xi_{3}|,$ $|\xi_{5}|\geq80|\xi_{6}|$
of  Lemma 3.6.

\noindent When $\frac{|\xi_{5}|}{80}\geq|\xi_{1}|\sim |\xi_{4}|\geq80 |\xi_{3}|,$ $|\xi_{5}|\leq 80|\xi_{6}|,$
this case can be proved similarly to case $\frac{|\xi_{5}|}{80}\geq|\xi_{1}|\geq80 |\xi_{2}|,$ $|\xi_{5}|\leq 80|\xi_{6}|$
of Lemma 3.6.

\noindent When $|\xi_{1}|\sim |\xi_{5}|\geq 80|\xi_{6}|,$ this case can be proved similarly to case
 $|\xi_{1}|\sim |\xi_{4}|,\frac{|\xi_{5}|}{80}\leq|\xi_{1}|\leq 80|\xi_{5}|$, $|\xi_{5}|\geq 80|\xi_{6}|$ of Lemma 3.6.

\noindent
When $|\xi_{1}|\sim |\xi_{6}|\geq 80|\xi_{7}|,$
this case can be proved similarly to case
 $|\xi_{1}|\sim |\xi_{4}|,\frac{|\xi_{5}|}{80}\leq|\xi_{1}|\leq 80|\xi_{5}|$, $|\xi_{5}|\geq 80|\xi_{6}|$
 of Lemma 3.6.

\noindent
When $|\xi_{1}|\sim |\xi_{7}|\geq 80|\xi_{8}|,$
this case can be proved similarly to case
 $|\xi_{1}|\sim |\xi_{4}|,\frac{|\xi_{5}|}{80}\leq|\xi_{1}|\leq 80|\xi_{5}|$, $|\xi_{5}|\geq 80|\xi_{6}|$
 of Lemma 3.6.

\noindent
When $|\xi_{1}|\sim |\xi_{8}|,$
this case can be proved similarly to case $|\xi_{1}|\sim |\xi_{4}|,\frac{|\xi_{5}|}{80}\leq|\xi_{1}|\leq 80|\xi_{5}|$,
$|\xi_{5}|\leq 80|\xi_{6}|,|\xi_{6}|\leq80|\xi_{7}|,|\xi_{7}|\leq 80|\xi_{8}|$ of Lemma 3.6.

This completes the proof of Lemma 3.6.

\begin{Lemma}\label{Lemma3.7}
Let $s\geq\frac{17}{112}+\epsilon$ and $\sigma=\frac{3}{14}+2\epsilon$ and
$z_{j}=\eta_{T}S(t)\phi^{\omega}(1\leq j\leq5,j\in N)$ and $v_{j}=\eta v(6\leq j\leq8,j\in N)$. Then, we have that
\begin{eqnarray*}
\left|\int_{\SR}\int_{\SR}J^{\sigma}\partial_{x}
\left(\left(\prod_{j=1}^{5}z_{j}\right)\left(\prod_{j=6}^{8}v_{j}\right)\right)h(x,t)dxdt\right|\leq CT^{-\frac{\epsilon}{50}}R^{5}\left(\prod_{j=6}^{8}\|v_{j}\|_{X_{\sigma,b}}\right)\|h\|_{X_{0,\frac{1}{2}-\frac{\epsilon}{12}}},
\end{eqnarray*}
outside a set of probability at most
$C{\rm exp}\left(-C^{'}\frac{R^{2}}{\|\phi\|_{H^{s}}^{2}}\right).
$
\end{Lemma}
\noindent {\bf Proof.}We dyadically decompose $z_{j}$ with $(1\leq j\leq5,j\in N)$ and
 $v_{j}$ with $(6\leq j\leq 8,j\in N)$ and $h$ such that frequency supports are
$\left\{|\xi_{j}|\sim N_{j}\right\}$ for some dyadically $N_{j}\geq 1$ and we still denote them by
$z_{j}$ with $(1\leq j\leq 5,j\in N)$ and $v_{j}$ with $(6\leq j\leq 8,j\in Z),$ and $h$.
 In this case, we divide the frequency  into
$$|\xi_{l}|\geq {\rm max}\left\{|\xi_{j}|,1\leq j\leq 5,j\neq l,j\in N\right\}(1\leq l\leq 5,l\in N).$$
Without loss of generality, we can assume that
$|\xi_{1}|\geq |\xi_{2}|\geq|\xi_{3}\geq |\xi_{5}|$ and $ |\xi_{6}|\geq|\xi_{7}\geq |\xi_{8}|$.

\noindent We define $$I_{7}=\left|\int_{\SR}\int_{\SR}J^{\sigma}\partial_{x}
\left(\left(\prod_{j=1}^{5}z_{j}\right)\left(\prod_{j=6}^{8}v_{j}\right)\right)h(x,t)dxdt\right|.$$

\noindent When $|\xi_{1}|\geq80 |\xi_{2}|,$ $|\xi_{1}|\geq 80|\xi_{6}|,$
by using the H\"older inequality, (\ref{2.04})-(\ref{2.05}), (\ref{2.016}), Lemmas 2.4,  2.1, 2.3,
   we get that
\begin{eqnarray*}
&&I_{7}\leq CN_{1}^{-\frac{1}{20}}
\left\|I^{1/2}I^{1/2}_{-}(J ^{s}z_{1},J ^{s}z_{2})\right\|_{L_{xt}^{2}}
\left(\prod_{j=4}^{5}\|z_{j}\|_{L_{xt}^{\infty}}\right)
\left(\prod_{j=6}^{8}\|J^{-\frac{2}{7}}v_{j}\|_{L_{xt}^{\infty}}\right)\nonumber\\&&\qquad \times
\|I^{\frac{1-\epsilon}{2}}I^{\frac{1-\epsilon}{2}}_{-}(J ^{s}z_{3},h)\|_{L_{xt}^{2}}\nonumber\\&&\leq
CN_{1}^{-\frac{1}{20}}\left(\prod_{j=1}^{3}\|z_{j}\|_{X_{s,c}}\right)
\left(\prod_{j=4}^{5}\|z_{j}\|_{L_{xt}^{\infty}}\right)\left(\prod_{j=6}^{8}\|v_{j}\|_{X_{\sigma,b}}\right)
\|h\|_{X_{0,\frac{1}{2}-\frac{\epsilon}{12}}}\nonumber\\
&&\leq CN_{1}^{-\frac{1}{20}}T^{-\frac{\epsilon}{100}}\left(\prod_{j=1}^{3}\|P_{N_{j}}\phi^{\omega}\|_{H^{s}}\right)
\left(\prod_{j=4}^{5}\|z_{j}\|_{L_{xt}^{\infty}}\right)\left(\prod_{j=6}^{8}\|v_{j}\|_{X_{\sigma,b}}\right)
\|h\|_{X_{0,\frac{1}{2}-\frac{\epsilon}{12}}}\nonumber\\&&\leq CT^{-\frac{\epsilon}{100}}R^{5}
\left(\prod_{j=6}^{8}\|v_{j}\|_{X_{\sigma,b}}\right)
\|h\|_{X_{0,\frac{1}{2}-\frac{\epsilon}{12}}},
\end{eqnarray*}
outside a set of probability at most
$
 C{\rm exp}\left(-C^{'}\frac{R^{2}}{\|\phi\|_{H^{s}}^{2}}\right).
$

\noindent When $|\xi_{1}|\geq80 |\xi_{2}|,$ $\frac{|\xi_{6}|}{80}\leq |\xi_{1}|\leq 80|\xi_{6}|,$
by using the H\"older inequality,
(\ref{2.05})-(\ref{2.07}), (\ref{2.014}),  Lemmas 2.4, 2.1, 2.2,
we have that
\begin{eqnarray*}
&&I_{7}\leq CN_{1}^{-\frac{1}{20}}
\left\|I^{1/2}I^{1/2}_{-}(J ^{s}z_{1},J ^{s}z_{2})\right\|_{L_{xt}^{2}}
\left(\prod_{j=3}^{5}\|z_{j}\|_{L_{xt}^{\frac{168}{6+21\epsilon}}}\right)
\|v_{6}\|_{L_{xt}^{\frac{28}{2-7\epsilon}}}\|J^{\sigma}v_{7}\|_{L_{xt}^{8}}\nonumber\\&&
\qquad\times\|v_{8}\|_{L_{xt}^{\frac{28}{2-7\epsilon}}}\|h\|_{L_{xt}^{\frac{8}{1+\epsilon}}}\nonumber\\&&\leq
CN_{1}^{-\frac{1}{20}}\left(\prod_{j=1}^{2}\|z_{j}\|_{X_{s,c}}\right)
\left(\prod_{j=3}^{5}\|z_{j}\|_{L_{xt}^{\frac{168}{6+21\epsilon}}}\right)\left(\prod_{j=6}^{8}\|v_{j}\|_{X_{\sigma,b}}\right)
\|h\|_{X_{0,\frac{1}{2}-\frac{\epsilon}{12}}}\nonumber\\
&&\leq CN_{1}^{-\frac{1}{20}}T^{-\frac{\epsilon}{50}}\left(\prod_{j=1}^{2}\|P_{N_{j}}\phi^{\omega}\|_{H^{s}}\right)
\left(\prod_{j=3}^{5}\|z_{j}\|_{L_{xt}^{\frac{168}{6+21\epsilon}}}\right)\left(\prod_{j=6}^{8}\|v_{j}\|_{X_{\sigma,b}}\right)
\|h\|_{X_{0,\frac{1}{2}-\frac{\epsilon}{12}}}\nonumber\\&&\leq C
T^{-\frac{\epsilon}{50}}R^{5}\left(\prod_{j=6}^{8}\|v_{j}\|_{X_{\sigma,b}}\right)
\|h\|_{X_{0,\frac{1}{2}-\frac{\epsilon}{12}}},
\end{eqnarray*}
outside a set of probability at most
$
C{\rm exp}\left(-C^{'}\frac{R^{2}}{\|\phi\|_{H^{s}}^{2}}\right).
$

\noindent When $\frac{|\xi_{6}|}{80}\geq|\xi_{1}|\geq80 |\xi_{2}|,$  by using the H\"older inequality,
(\ref{2.05})-(\ref{2.07}), Lemmas 2.4, 2.1, 2.2,
we have that
\begin{eqnarray*}
&&I_{7}\leq C
\left\|I^{1/2}I^{1/2}_{-}(J ^{\sigma}v_{6},z_{1})\right\|_{L_{xt}^{2}}
\left(\prod_{j=2}^{5}\|z_{j}\|_{L_{xt}^{\frac{224}{13+21\epsilon}}}\right)
\left(\prod_{j=7}^{8}\|v_{j}\|_{L_{xt}^{\frac{28}{2-7\epsilon}}}\right)
\|h\|_{L_{xt}^{\frac{8}{1+\epsilon}}}\nonumber\\&&\leq
CN_{1}^{-2\epsilon}\|z_{1}\|_{X_{0,c}}\left(\prod_{j=2}^{5}\|z_{j}\|_{L_{xt}^{\frac{224}{13+21\epsilon}}}\right)
\left(\prod_{j=6}^{8}\|v_{j}\|_{X_{\sigma,b}}\right)
\|h\|_{X_{0,\frac{1}{2}-\frac{\epsilon}{12}}}\nonumber\\
&&\leq CT^{-\frac{\epsilon}{100}}N_{1}^{-2\epsilon}\|P_{N_{1}}\phi^{\omega}\|_{L^{2}}
\left(\prod_{j=2}^{5}\|z_{j}\|_{L_{xt}^{\frac{224}{13+21\epsilon}}}\right)\left(\prod_{j=6}^{8}\|v_{j}\|_{X_{\sigma,b}}\right)
\|h\|_{X_{0,\frac{1}{2}-\frac{\epsilon}{12}}}\nonumber\\
&&\leq CT^{-\frac{\epsilon}{100}}R^{5}\left(\prod_{j=6}^{8}\|v_{j}\|_{X_{\sigma,b}}\right)
\|h\|_{X_{0,\frac{1}{2}-\frac{\epsilon}{12}}},
\end{eqnarray*}
outside a set of probability at most
$
 C{\rm exp}\left(-C^{'}\frac{R^{2}}{\|\phi\|_{H^{\epsilon}}^{2}}\right).
$

\noindent When $|\xi_{1}|\sim |\xi_{2}|\geq 80|\xi_{3}|,$ $|\xi_{1}|\geq \frac{|\xi_{6}|}{80},$ by using the H\"older inequality,
(\ref{2.05})-(\ref{2.07}), Lemmas 2.4, 2.1, 2.2,
we have that
\begin{eqnarray*}
&&I_{7}\leq CN_{2}^{-\frac{1}{20}}
\|J^{s}z_{1}\|_{L_{xt}^{\frac{168}{9+35\epsilon}}}\left\|I^{1/2}I^{1/2}_{-}(J^{s}z_{2},J^{s}z_{3})\right\|_{L_{xt}^{2}}
\left(\prod_{j=4}^{5}\|J^{s}z_{j}\|_{L_{xt}^{\frac{168}{9+35\epsilon}}}\right)\nonumber\\&&\qquad\times\left(\prod_{j=6}^{8}
\|v_{j}\|_{L_{xt}^{\frac{28}{2-7\epsilon}}}\right)
\|h\|_{L_{xt}^{\frac{8}{1+\epsilon}}}\nonumber\\&&\leq
CN_{2}^{-\frac{1}{20}}\|J^{s}z_{1}\|_{L_{xt}^{\frac{168}{9+35\epsilon}}}\left(\prod_{j=2}^{3}\|z_{j}\|_{X_{s,c}}\right)
\left(\prod_{j=4}^{5}\|J^{s}z_{j}\|_{L_{xt}^{\frac{168}{9+35\epsilon}}}\right)\left(\prod_{j=6}^{8}\|v_{j}\|_{X_{\sigma,b}}\right)
\|h\|_{X_{0,\frac{1}{2}-\frac{\epsilon}{12}}}\nonumber\\
&&\leq CN_{2}^{-\frac{1}{20}}T^{-\frac{\epsilon}{50}}\|J^{s}z_{1}\|_{L_{xt}^{\frac{168}{9+35\epsilon}}}
\left(\prod_{j=2}^{3}\|P_{N_{j}}\phi^{\omega}\|_{H^{s}}\right)
\left(\prod_{j=4}^{5}\|J^{s}z_{j}\|_{L_{xt}^{\frac{168}{9+35\epsilon}}}\right)\nonumber\\&&\qquad\times
\left(\prod_{j=6}^{8}\|v_{j}\|_{X_{\sigma,b}}\right)
\|h\|_{X_{0,\frac{1}{2}-\frac{\epsilon}{12}}}\leq CT^{-\frac{\epsilon}{50}}R^{5}\left(\prod_{j=6}^{8}\|v_{j}\|_{X_{\sigma,b}}\right)
\|h\|_{X_{0,\frac{1}{2}-\frac{\epsilon}{12}}},
\end{eqnarray*}
outside a set of probability at most
$
 C{\rm exp}\left(-C^{'}\frac{R^{2}}{\|\phi\|_{H^{s}}^{2}}\right).
$

\noindent When $\frac{|\xi_{6}|}{80}\geq|\xi_{1}|\sim |\xi_{2}|\geq 80|\xi_{3}|,$  this case can be proved similarly to case
$\frac{|\xi_{6}|}{80}\geq|\xi_{1}|\geq 80|\xi_{2}|$ of Lemma 3.7.

\noindent When $|\xi_{1}|\sim |\xi_{3}|\geq 80|\xi_{4}|,$ $|\xi_{1}|\geq \frac{|\xi_{6}|}{80},$
this case can be proved similarly to case
$|\xi_{1}|\sim |\xi_{2}|\geq 80|\xi_{3}|,$ $|\xi_{1}|\geq \frac{|\xi_{6}|}{80}$ of Lemma 3.7.

\noindent When $\frac{|\xi_{6}|}{80}\geq|\xi_{1}|\sim |\xi_{3}|\geq 80|\xi_{4}|,$  this case can be proved similarly to case
$\frac{|\xi_{6}|}{80}\geq|\xi_{1}|\geq 80|\xi_{2}|$ of Lemma 3.7.

\noindent When $|\xi_{1}|\sim |\xi_{4}|\geq 80|\xi_{5}|,$ $|\xi_{1}|\geq \frac{|\xi_{6}|}{80}$,
this case can be proved similarly to case
$|\xi_{1}|\sim |\xi_{2}|\geq 80|\xi_{3}|,$ $|\xi_{1}|\geq \frac{|\xi_{6}|}{80}$ of Lemma 3.7.

\noindent When $\frac{|\xi_{6}|}{80}\geq|\xi_{1}|\sim |\xi_{4}|\geq 80|\xi_{5}|,$
this case can be proved similarly to case $\frac{|\xi_{6}|}{80}\geq|\xi_{1}|\geq 80|\xi_{2}|$ of Lemma 3.7.

\noindent When $|\xi_{1}|\sim |\xi_{5}|,|\xi_{1}|\geq 80|\xi_{6}|$,  by using the H\"older inequality,
(\ref{2.05})-(\ref{2.07}), Lemmas 2.4, 2.1, 2.2,  we have that
\begin{eqnarray*}
&&I_{7}
\leq CN_{5}^{-\frac{1}{20}}
\left(\prod_{j=2}^{5}\|J^{s}z_{j}\|_{L_{xt}^{\frac{224}{13+21\epsilon}}}\right)
\left\|I^{1/2}I^{1/2}_{-}(J^{\sigma}v_{6},J^{s}z_{1})\right\|_{L_{xt}^{2}}
\left(\prod_{j=7}^{8}\|v_{j}\|_{L_{xt}^{\frac{28}{2-7\epsilon}}}\right)
\|h\|_{L_{xt}^{\frac{8}{1+\epsilon}}}\nonumber\\&&\leq CN_{5}^{-\frac{1}{20}}\|z_{1}\|_{X_{s,c}}
\left(\prod_{j=2}^{5}\|J^{s}z_{j}\|_{L_{xt}^{\frac{224}{13+21\epsilon}}}\right)
\left(\prod_{j=6}^{8}\|v_{j}\|_{X_{\sigma,b}}\right)\|h\|_{X_{0,\frac{1}{2}-\frac{\epsilon}{12}}}\nonumber\\
&&\leq CN_{5}^{-\frac{1}{20}}T^{-\frac{\epsilon}{100}}\|P_{N_{1}}\phi^{\omega}\|_{H^{s}}\left(\prod_{j=2}^{5}
\|J^{s}z_{j}\|_{L_{xt}^{\frac{224}{13+21\epsilon}}}\right)
\left(\prod_{j=6}^{8}\|v_{j}\|_{X_{\sigma,b}}\right)\|h\|_{X_{0,\frac{1}{2}-\frac{\epsilon}{12}}}\nonumber\\
&&\leq CT^{-\frac{\epsilon}{100}}R^{5}
\left(\prod_{j=6}^{8}\|v_{j}\|_{X_{\sigma,b}}\right)\|h\|_{X_{0,\frac{1}{2}-\frac{\epsilon}{12}}},
\end{eqnarray*}
outside a set of probability at most
$
C{\rm exp}\left(-C^{'}\frac{R^{2}}{\|\phi\|_{H^{s}}^{2}}\right).
$

\noindent When $|\xi_{1}|\sim |\xi_{5}|,\frac{|\xi_{6}|}{80}\leq |\xi_{1}|\leq 80|\xi_{6}|$,
$|\xi_{6}|\geq 80|\xi_{7}|,$ by using (\ref{2.05})-(\ref{2.07}),  Lemma 2.2,  we have that
\begin{eqnarray*}
&&I_{7}
\leq CN_{5}^{-\frac{1}{20}}
\left(\prod_{j=1}^{5}\|J^{s}z_{j}\|_{L_{xt}^{\frac{280}{17+7\epsilon}}}\right)
\left\|I^{1/2}I^{1/2}_{-}(J^{\sigma}v_{6},J^{s}v_{7})\right\|_{L_{xt}^{2}}
\|v_{8}\|_{L_{xt}^{\frac{28}{2-7\epsilon}}}
\|h\|_{L_{xt}^{\frac{8}{1+\epsilon}}}\nonumber\\&&\leq CN_{5}^{-\frac{1}{20}}
\left(\prod_{j=1}^{5}\|J^{s}z_{j}\|_{L_{xt}^{\frac{280}{17+7\epsilon}}}\right)
\left(\prod_{j=6}^{8}\|v_{j}\|_{X_{\sigma,b}}\right)
\|h\|_{X_{0,\frac{1}{2}-\frac{\epsilon}{12}}}\nonumber\\
&&\leq CN_{5}^{-\frac{1}{20}}T^{-\frac{\epsilon}{100}}\left(\prod_{j=1}^{5}
\|J^{s}z_{j}\|_{L_{xt}^{\frac{280}{17+7\epsilon}}}\right)
\left(\prod_{j=6}^{8}\|v_{j}\|_{X_{\sigma,b}}\right)
\|h\|_{X_{0,\frac{1}{2}-\frac{\epsilon}{12}}}\nonumber\\
&&\leq CT^{-\frac{\epsilon}{100}}R^{5}
\left(\prod_{j=6}^{8}\|v_{j}\|_{X_{\sigma,b}}\right)
\|h\|_{X_{0,\frac{1}{2}-\frac{\epsilon}{12}}},
\end{eqnarray*}
outside a set of probability at most
$
 C{\rm exp}\left(-C^{'}\frac{R^{2}}{\|\phi\|_{H^{\epsilon}}^{2}}\right).
$

\noindent When $|\xi_{1}|\sim |\xi_{5}|,\frac{|\xi_{6}|}{80}\leq |\xi_{5}|\leq 80|\xi_{6}|$,
$|\xi_{6}|\leq 80|\xi_{7}|,$ $|\xi_{7}|\geq 80|\xi_{8}|,$ this case can be proved similarly to case
$|\xi_{1}|\sim |\xi_{5}|,\frac{|\xi_{6}|}{80}\leq |\xi_{5}|\leq 80|\xi_{6}|$,
$|\xi_{6}|\geq 80|\xi_{7}|$ of Lemma 3.7.

\noindent When $|\xi_{1}|\sim |\xi_{5}|,\frac{|\xi_{6}|}{80}\leq |\xi_{1}|\leq 80|\xi_{6}|$,
 $|\xi_{6}|\leq 80|\xi_{7}|,$ $|\xi_{7}|\leq 80|\xi_{8}|,$
by using the H\"older  inequality,  (\ref{2.06}),
(\ref{2.014}),  Lemmas 2.2,   we get that
\begin{eqnarray*}
&&I_{7}\leq CN_{1}^{-\frac{1}{20}}
\left(\prod_{j=1}^{5}\|J^{s}z_{j}\|_{L_{xt}^{\frac{40}{4-\epsilon}}}\right)\left(\prod_{j=6}^{8}\|J^{\sigma}
v_{j}\|_{L_{xt}^{8}}\right)
\|h\|_{L_{xt}^{\frac{8}{1+\epsilon}}}\nonumber\\&&\leq  CN_{1}^{-\frac{1}{20}}\left(\prod_{j=1}^{5}\|J^{s}z_{j}\|_{L_{xt}^{\frac{40}{4-\epsilon}}}\right)\left(\prod_{j=6}^{8}
\|v_{j}\|_{X_{\sigma,b}}\right)\|h\|_{X_{0,\frac{1}{2}-\frac{\epsilon}{12}}}\nonumber\\
&&\leq CR^{5}\left(\prod_{j=6}^{8}
\|v_{j}\|_{X_{\sigma,b}}\right)\|h\|_{X_{0,\frac{1}{2}-\frac{\epsilon}{12}}},
\end{eqnarray*}
outside a set of probability at most
$
C{\rm exp}\left(-C^{'}\frac{R^{2}}{T^{\frac{4-\epsilon}{20}}\|\phi\|_{H^{s}}^{2}}\right).
$

\noindent When $\frac{|\xi_{6}|}{80}\geq|\xi_{1}|\sim |\xi_{5}|$,
this case can be proved similarly to case  $\frac{|\xi_{6}|}{80}\geq|\xi_{1}|\geq80|\xi_{2}|$ of Lemma 3.7.

\noindent When $|\xi_{1}|\sim |\xi_{6}|, |\xi_{1}|\geq80|\xi_{7}|$, this case can
be proved similarly to case $|\xi_{1}|\sim |\xi_{5}|,|\xi_{1}|\geq 80|\xi_{6}|$ of Lemma 3.7.

\noindent When $|\xi_{1}|\sim |\xi_{6}|,\frac{|\xi_{7}|}{80} \leq |\xi_{1}|\leq80|\xi_{7}|$,
  $|\xi_{7}|\geq 80|\xi_{8}|,$
this case can be proved similarly to case
$|\xi_{1}|\sim |\xi_{5}|,\frac{|\xi_{6}|}{80}\leq |\xi_{5}|\leq 80|\xi_{6}|$,
$|\xi_{6}|\geq 80|\xi_{7}|$ of Lemma 3.7.

\noindent When $\frac{|\xi_{7}|}{80}\geq|\xi_{1}|\sim |\xi_{6}|$,
  this case can be proved similarly to case  $\frac{|\xi_{6}|}{80}\geq|\xi_{1}|\geq80|\xi_{2}|$ of  Lemma 3.7.

\noindent When $|\xi_{1}|\sim |\xi_{7}|, |\xi_{1}|\geq80|\xi_{8}|$,   this case can be proved similarly to
 case $|\xi_{1}|\sim |\xi_{5}|,|\xi_{1}|\geq 80|\xi_{6}|$ of Lemma 3.7.

\noindent When $|\xi_{1}|\sim |\xi_{7}|, \frac{|\xi_{8}|}{80}\leq|\xi_{1}|\leq80|\xi_{8}|$,
 this case can be proved similarly to case  $|\xi_{1}|\sim |\xi_{5}|,\frac{|\xi_{6}|}{80}\leq |\xi_{1}|\leq 80|\xi_{6}|$,
 $|\xi_{6}|\leq 80|\xi_{7}|,$ $|\xi_{7}|\leq 80|\xi_{8}|$ of Lemma 3.7.

\noindent When $\frac{|\xi_{8}|}{80}\geq|\xi_{1}|\sim |\xi_{7}|$,
 this case can be proved similarly to $\frac{|\xi_{7}|}{80}\geq|\xi_{1}|\sim |\xi_{6}|$  of  Lemma 3.7.

\noindent When $|\xi_{1}|\sim|\xi_{8}|$, this case can be proved similarly to case
$|\xi_{1}|\sim |\xi_{5}|,\frac{|\xi_{6}|}{80}\leq |\xi_{5}|\leq 80|\xi_{6}|$,
 $|\xi_{6}|\leq 80|\xi_{7}|,$ $|\xi_{7}|\leq 80|\xi_{8}|$ of  Lemma 3.7.

This completes the proof of Lemma 3.7.

\begin{Lemma}\label{Lemma3.8}
Let $s\geq\frac{17}{112}+\epsilon$ and $\sigma=\frac{3}{14}+2\epsilon$
and $z_{j}=\eta_{T}S(t)\phi^{\omega}(1\leq j\leq6,j\in N)$ and $v_{j}=\eta v$ with $j=7,8$. Then, we have that
\begin{eqnarray*}
\left|\int_{\SR}\int_{\SR}J^{\sigma}\partial_{x}\left(\left(\prod_{j=1}^{6}z_{j}\right)
\left(\prod_{j=7}^{8}v_{j}\right)\right)h(x,t)dxdt\right|\leq CT^{-\frac{3\epsilon}{100}}R^{6}\left(\prod_{j=7}^{8}\|v_{j}\|_{X_{\sigma,b}}\right)
\|h\|_{X_{0,\frac{1}{2}-\frac{\epsilon}{12}}},
\end{eqnarray*}
outside a set of probability at most
$
C{\rm exp}\left(-C^{'}\frac{R^{2}}{\|\phi\|_{H^{s}}^{2}}\right).
$
\end{Lemma}
\noindent {\bf Proof.} We dyadically decompose $z_{j}$ with $(1\leq j\leq6,j\in N)$
and $v_{j}$ with $(j=7,8),$ and $h$ such that frequency supports are
$\left\{|\xi_{j}|\sim N_{j}\right\}$ for some dyadically $N_{j}\geq 1$ and we still
denote them by $z_{j}$ with $(1\leq j\leq 6,j\in N)$ and $h$.
In this case, we divide the frequency  into
$$|\xi_{l}|\geq {\rm max}\left\{|\xi_{j}|,1\leq j\leq 6,j\neq l,j\in N\right\}(1\leq l\leq 6,l\in N).$$
Without loss of generality, we can assume that $|\xi_{1}|\geq |\xi_{2}|
\geq|\xi_{3}\geq |\xi_{5}|\geq|\xi_{6}|$ and $ |\xi_{7}\geq |\xi_{8}|$.

\noindent We define $$I_{8}=\left|\int_{\SR}\int_{\SR}J^{\sigma}\partial_{x}\left(\left(\prod_{j=1}^{6}z_{j}\right)
\left(\prod_{j=7}^{8}v_{j}\right)\right)h(x,t)dxdt\right|.$$

\noindent When $|\xi_{1}|\geq80 |\xi_{2}|,$ $|\xi_{1}|\geq 80|\xi_{7}|,$
by using the H\"older inequality,  (\ref{2.04})-(\ref{2.05}), (\ref{2.016}), Lemmas 2.4, 2.1, 2.3,
we have that
\begin{eqnarray*}
&&I_{8}\leq CN_{1}^{-\frac{1}{10}}
\left\|I^{1/2}I^{1/2}_{-}(J ^{s}z_{1},J ^{s}z_{2})\right\|_{L_{xt}^{2}}\left(\prod_{j=4}^{6}\|z_{j}\|_{L_{xt}^{\infty}}\right)
\left(\prod_{j=7}^{8}\|J^{-\frac{2}{7}}v_{j}\|_{L_{xt}^{\infty}}\right)\nonumber\\&&\qquad\times
\|I^{\frac{1-\epsilon}{2}}I^{\frac{1-\epsilon}{2}}_{-}(J ^{s}z_{3},h)\|_{L_{xt}^{2}}\nonumber\\&&\leq
CN_{1}^{-\frac{1}{10}}\left(\prod_{j=1}^{3}\|z_{j}\|_{X_{s,c}}\right)
\left(\prod_{j=4}^{6}\|z_{j}\|_{L_{xt}^{\infty}}\right)\left(\prod_{j=7}^{8}\|v_{j}\|_{X_{\sigma,b}}\right)
\|h\|_{X_{0,\frac{1}{2}-\frac{\epsilon}{12}}}\nonumber\\
&&\leq CN_{1}^{-\frac{1}{10}}T^{-\frac{3\epsilon}{100}}
\left(\prod_{j=1}^{3}\|P_{N_{j}}\phi^{\omega}\|_{H^{s}}\right)\left(\prod_{j=4}^{6}\|z_{j}\|_{L_{xt}^{\infty}}\right)
\left(\prod_{j=7}^{8}\|v_{j}\|_{X_{\sigma,b}}\right)
\|h\|_{X_{0,\frac{1}{2}-\frac{\epsilon}{12}}}\nonumber\\
&&\leq CT^{-\frac{3\epsilon}{100}}R^{6}\left(\prod_{j=7}^{8}\|v_{j}\|_{X_{\sigma,b}}\right)
\|h\|_{X_{0,\frac{1}{2}-\frac{\epsilon}{12}}},
\end{eqnarray*}
outside a set of probability at most
$
C{\rm exp}\left(-C^{'}\frac{R^{2}}{\|\phi\|_{H^{s}}^{2}}\right).
$

\noindent When $|\xi_{1}|\geq80 |\xi_{2}|,$ $\frac{|\xi_{7}|}{80}\leq|\xi_{1}|\leq 80|\xi_{7}|,$
by using the H\"older inequality,    (\ref{2.05})-(\ref{2.06}), (\ref{2.014}), Lemmas 2.4, 2.1, 2.2,
we have that
\begin{eqnarray*}
&&I_{8}\leq CN_{1}^{-\frac{1}{20}}
\left\|I^{1/2}I^{1/2}_{-}(J ^{s}z_{1},J ^{s}z_{2})\right\|_{L_{xt}^{2}}
\left(\prod_{j=3}^{6}\|z_{j}\|_{L_{xt}^{\frac{32}{1-\epsilon}}}\right)
\left(\prod_{j=7}^{8}\|J^{\sigma}v_{j}\|_{L_{xt}^{8}}\right)
\|h\|_{L_{xt}^{\frac{8}{1+\epsilon}}}\nonumber\\&&\leq
CN_{1}^{-\frac{1}{20}}\left(\prod_{j=1}^{2}\|z_{j}\|_{X_{s,c}}\right)\left(\prod_{j=3}^{6}
\|z_{j}\|_{L_{xt}^{\frac{32}{1-\epsilon}}}\right)\left(\prod_{j=7}^{8}\|v_{j}\|_{X_{\sigma,b}}\right)
\|h\|_{X_{0,\frac{1}{2}-\frac{\epsilon}{12}}}\nonumber\\
&&\leq CN_{1}^{-\frac{1}{20}}T^{-\frac{\epsilon}{50}}
\left(\prod_{j=1}^{2}\|P_{N_{j}}\phi^{\omega}\|_{H^{s}}\right)\left(\prod_{j=3}^{6}
\|z_{j}\|_{L_{xt}^{\frac{224}{10+7\epsilon}}}\right)\left(\prod_{j=7}^{8}\|v_{j}\|_{X_{\sigma,b}}\right)
\|h\|_{X_{0,\frac{1}{2}-\frac{\epsilon}{12}}}\nonumber\\
&&\leq CT^{-\frac{\epsilon}{50}}R^{6}\left(\prod_{j=7}^{8}\|v_{j}\|_{X_{\sigma,b}}\right)
\|h\|_{X_{0,\frac{1}{2}-\frac{\epsilon}{12}}},
\end{eqnarray*}
outside a set of probability at most
$ C{\rm exp}\left(-C^{'}\frac{R^{2}}{\|\phi\|_{H^{s}}^{2}}\right).
$

\noindent When $\frac{|\xi_{7}|}{80}\geq|\xi_{1}|\geq80 |\xi_{2}|,$  by using the  H\"older  inequality,
   (\ref{2.05})-(\ref{2.07}), Lemmas 2.4, 2.1, 2.2,
we have that
\begin{eqnarray*}
&&I_{8}\leq C
\left\|I^{1/2}I^{1/2}_{-}(J ^{\sigma}v_{7},z_{1})\right\|_{L_{xt}^{2}}
\left(\prod_{j=2}^{6}\|z_{j}\|_{L_{xt}^{\frac{280}{17+7\epsilon}}}\right)
\|v_{8}\|_{L_{xt}^{\frac{28}{2-7\epsilon}}}
\|h\|_{L_{xt}^{\frac{8}{1+\epsilon}}}\nonumber\\&&\leq
C\|z_{1}\|_{X_{0,c}}
\left(\prod_{j=2}^{6}\|z_{j}\|_{L_{xt}^{\frac{280}{17+7\epsilon}}}\right)
\left(\prod_{j=7}^{8}\|v_{j}\|_{X_{\sigma,b}}\right)
\|h\|_{X_{0,\frac{1}{2}-\frac{\epsilon}{12}}}\nonumber\\
&&\leq CT^{-\frac{\epsilon}{100}}\|P_{N_{1}}\phi^{\omega}\|_{L^{2}}
\left(\prod_{j=2}^{6}\|z_{j}\|_{L_{xt}^{\frac{280}{17+7\epsilon}}}\right)
\left(\prod_{j=7}^{8}\|v_{j}\|_{X_{\sigma,b}}\right)
\|h\|_{X_{0,\frac{1}{2}-\frac{\epsilon}{12}}}\nonumber\\
&&\leq CT^{-\frac{\epsilon}{100}}R^{6}\left(\prod_{j=7}^{8}\|v_{j}\|_{X_{\sigma,b}}\right)
\|h\|_{X_{0,\frac{1}{2}-\frac{\epsilon}{12}}},
\end{eqnarray*}
outside a set of probability at most
$
C{\rm exp}\left(-C^{'}\frac{R^{2}}{\|\phi\|_{H^{\epsilon}}^{2}}\right).
$

\noindent When $|\xi_{1}|\sim |\xi_{2}|\geq 80|\xi_{3}|,$
$|\xi_{1}|\geq \frac{|\xi_{7}|}{80},$ by using the H\"older inequality,
 (\ref{2.05})-(\ref{2.06}), (\ref{2.014}),  Lemmas 2.4, 2.1, 2.2,
we obtain that
\begin{eqnarray*}
&&I_{8}
\leq CN_{1}^{-\frac{1}{20}}\|J^{s}z_{1}\|_{L_{xt}^{\frac{32}{1-\epsilon}}}
\left\|I^{1/2}I^{1/2}_{-}(J^{s}z_{2},J^{s}z_{3})\right\|_{L_{xt}^{2}}
\left(\prod_{j=4}^{6}\|J^{s}z_{j}\|_{L_{xt}^{\frac{32}{1-\epsilon}}}\right)\left(\prod_{j=7}^{8}
\|J^{\sigma}v_{j}\|_{L_{xt}^{8}}\right)\nonumber\\&&\qquad \times
\|h\|_{L_{xt}^{\frac{8}{1+\epsilon}}}\nonumber\\&&\leq
CN_{1}^{-\frac{1}{20}}\|J^{s}z_{1}\|_{L_{xt}^{\frac{32}{1-\epsilon}}}
\left(\prod_{j=2}^{3}\|z_{j}\|_{X_{s,c}}\right)
\left(\prod_{j=4}^{6}\|J^{s}z_{j}\|_{L_{xt}^{\frac{32}{1-\epsilon}}}\right)
\left(\prod_{j=7}^{8}\|v_{j}\|_{X_{\sigma,b}}\right)
\|h\|_{X_{0,\frac{1}{2}-\frac{\epsilon}{12}}}\nonumber\\
&&\leq CN_{1}^{-\frac{1}{20}}T^{-\frac{\epsilon}{50}}\|J^{s}z_{1}\|_{L_{xt}^{\frac{32}{1-\epsilon}}}
\left(\prod_{j=2}^{3}\|P_{N_{j}}\phi^{\omega}\|_{H^{s}}\right)
\left(\prod_{j=4}^{6}\|J^{s}z_{j}\|_{L_{xt}^{\frac{32}{1-\epsilon}}}\right)
\nonumber\\&&\qquad\times\left(\prod_{j=7}^{8}\|v_{j}\|_{X_{\sigma,b}}\right)
\|h\|_{X_{0,\frac{1}{2}-\frac{\epsilon}{12}}}\leq CT^{-\frac{\epsilon}{50}}
R^{6}\left(\prod_{j=7}^{8}\|v_{j}\|_{X_{\sigma,b}}\right)
\|h\|_{X_{0,\frac{1}{2}-\frac{\epsilon}{12}}},
\end{eqnarray*}
outside a set of probability at most
$
 C{\rm exp}\left(-C^{'}\frac{R^{2}}{\|\phi\|_{H^{s}}^{2}}\right).
$

\noindent When $\frac{|\xi_{7}|}{80}\geq|\xi_{1}|\sim |\xi_{2}|\geq 80|\xi_{3}|,$
by using the  H\"older inequality,
 (\ref{2.05})-(\ref{2.07}), Lemmas 2.4, 2.1, 2.2,
we have that
\begin{eqnarray*}
&&I_{8}\leq C
\left\|I^{1/2}I^{1/2}_{-}(J^{\sigma}v_{7},z_{1})\right\|_{L_{xt}^{2}}
\left(\prod_{j=2}^{6}\|z_{j}\|_{L_{xt}^{\frac{280}{17+7\epsilon}}}\right)
\|v_{8}\|_{L_{xt}^{\frac{28}{2-7\epsilon}}}
\|h\|_{L_{xt}^{\frac{8}{1+\epsilon}}}\nonumber\\&&\leq
C\|z_{1}\|_{X_{0,c}}
\left(\prod_{j=2}^{6}\|z_{j}\|_{L_{xt}^{\frac{280}{17+7\epsilon}}}\right)
\left(\prod_{j=7}^{8}\|v_{j}\|_{X_{\sigma,b}}\right)
\|h\|_{X_{0,\frac{1}{2}-\frac{\epsilon}{12}}}\nonumber\\
&&\leq CT^{-\frac{\epsilon}{100}}\|P_{N_{1}}\phi^{\omega}\|_{L^{2}}
\left(\prod_{j=2}^{6}\|z_{j}\|_{L_{xt}^{\frac{280}{17+7\epsilon}}}\right)
\left(\prod_{j=7}^{8}\|v_{j}\|_{X_{\sigma,b}}\right)
\|h\|_{X_{0,\frac{1}{2}-\frac{\epsilon}{12}}}\nonumber\\&&\leq CT^{-\frac{\epsilon}{100}}R^{6}
\left(\prod_{j=7}^{8}\|v_{j}\|_{X_{\sigma,b}}\right)
\|h\|_{X_{0,\frac{1}{2}-\frac{\epsilon}{12}}},
\end{eqnarray*}
outside a set of probability at most
$
C{\rm exp}\left(-C^{'}\frac{R^{2}}{\|\phi\|_{H^{\epsilon}}^{2}}\right).
$

\noindent When $|\xi_{1}|\sim |\xi_{3}|\geq 80|\xi_{4}|,$ $|\xi_{1}|\geq \frac{|\xi_{7}|}{80},$
this case can be proved similarly to case  $|\xi_{1}|\sim |\xi_{2}|\geq 80|\xi_{3}|,$
$|\xi_{1}|\geq \frac{|\xi_{7}|}{80}$ of Lemma 3.8.

\noindent When $\frac{|\xi_{7}|}{80}\geq|\xi_{1}|\sim |\xi_{3}|\geq 80|\xi_{4}|,$
this case can be proved similarly to  case $\frac{|\xi_{7}|}{80}\geq|\xi_{1}|
\sim |\xi_{2}|\geq 80|\xi_{3}|$ of Lemma 3.8.

\noindent When $|\xi_{1}|\sim |\xi_{4}|\geq 80|\xi_{5}|,$ $|\xi_{1}|\geq \frac{|\xi_{7}|}{80}$,
this case can be proved similarly to case  $|\xi_{1}|\sim |\xi_{2}|\geq 80|\xi_{3}|,$
$|\xi_{1}|\geq \frac{|\xi_{7}|}{80}$ of Lemma 3.8.

\noindent When $\frac{|\xi_{7}|}{80}\geq|\xi_{1}|\sim |\xi_{4}|\geq 80|\xi_{5}|,$
this case can be proved similarly to  case $\frac{|\xi_{7}|}{80}\geq|\xi_{1}|\sim |\xi_{2}|
\geq 80|\xi_{3}|$ of Lemma 3.8.

\noindent When $|\xi_{1}|\sim |\xi_{5}|\geq80|\xi_{6}|,|\xi_{1}|\geq\frac{|\xi_{7}|}{80}$,
 this case can be proved similarly to case  $|\xi_{1}|\sim |\xi_{2}|\geq 80|\xi_{3}|,$
$|\xi_{1}|\geq \frac{|\xi_{7}|}{80}$ of Lemma 3.8.

\noindent When $\frac{|\xi_{7}|}{80}\geq|\xi_{1}|\sim |\xi_{5}|\geq80|\xi_{6}|$,
this case can be proved similarly to  case $\frac{|\xi_{7}|}{80}\geq|\xi_{1}|\sim |\xi_{2}|
\geq 80|\xi_{3}|$ of Lemma 3.8.

\noindent When $|\xi_{1}|\sim |\xi_{6}|, |\xi_{1}|\geq80|\xi_{7}|$,
by using the H\"older inequality, (\ref{2.05})-(\ref{2.07}),
Lemmas 2.4, 2.1, 2.2,  we have that
\begin{eqnarray*}
&&I_{8}
\leq CN_{6}^{-\frac{1}{20}}\left(\prod_{j=2}^{6}\|J^{s}z_{j}\|_{L_{xt}^{\frac{280}{17+7\epsilon}}}\right)
\left\|I^{1/2}I^{1/2}_{-}(J^{s}z_{1},J^{\sigma}v_{7})\right\|_{L_{xt}^{2}}
\|v_{8}\|_{L_{xt}^{\frac{28}{2-7\epsilon}}}
\|h\|_{L_{xt}^{\frac{8}{1+\epsilon}}}\nonumber\\&&\leq CN_{6}^{-\frac{1}{20}}
\|z_{1}\|_{X_{s,c}}\left(\prod_{j=2}^{6}\|J^{s}z_{j}\|_{L_{xt}^{\frac{280}{17+7\epsilon}}}\right)
\left(\prod_{j=7}^{8}\|v_{j}\|_{X_{\sigma,b}}\right)
\|h\|_{X_{0,\frac{1}{2}-\frac{\epsilon}{12}}}\nonumber\\
&&\leq CT^{-\frac{\epsilon}{100}}N_{6}^{-\frac{1}{20}}\|P_{N_{1}}\phi^{\omega}\|_{H^{s}}
\left(\prod_{j=2}^{6}\|J^{s}z_{j}\|_{L_{xt}^{\frac{280}{17+7\epsilon}}}\right)
\left(\prod_{j=7}^{8}\|v_{j}\|_{X_{\sigma,b}}\right)
\|h\|_{X_{0,\frac{1}{2}-\frac{\epsilon}{12}}}\nonumber\\
&&\leq CT^{-\frac{\epsilon}{100}}R^{6}
\left(\prod_{j=7}^{8}\|v_{j}\|_{X_{\sigma,b}}\right),
\end{eqnarray*}
outside a set of probability  at most
$
C{\rm exp}\left(-C^{'}\frac{R^{2}}{\|\phi\|_{H^{s}}^{2}}\right).
$

\noindent When $|\xi_{1}|\sim |\xi_{6}|, \frac{|\xi_{7}|}{80}
\leq |\xi_{1}|\leq80|\xi_{7}|$, $|\xi_{7}|\geq 80|\xi_{8}|,$
by using the H\"older inequality, (\ref{2.05})-(\ref{2.06}),
Lemma 2.2,  we conclude that
\begin{eqnarray*}
&&I_{8}
\leq CN_{6}^{-\frac{1}{20}}\left(\prod_{j=1}^{6}\|J^{s}z_{j}\|_{L_{xt}^{\frac{48}{3-\epsilon}}}\right)
\left\|I^{1/2}I^{1/2}_{-}(J^{\sigma}v_{7},J^{\sigma}v_{8})\right\|_{L_{xt}^{2}}
\|h\|_{L_{xt}^{\frac{8}{1+\epsilon}}}\nonumber\\&&\leq CN_{6}^{-\frac{1}{20}}
\left(\prod_{j=1}^{6}\|J^{s}z_{j}\|_{L_{xt}^{\frac{48}{3-\epsilon}}}\right)
\left(\prod_{j=7}^{8}\|v_{j}\|_{X_{\sigma,b}}\right)\|h\|_{X_{0,\frac{1}{2}-\frac{\epsilon}{12}}}\nonumber\\
&&\leq CR^{6}
\left(\prod_{j=7}^{8}\|v_{j}\|_{X_{\sigma,b}}\right)\|h\|_{X_{0,\frac{1}{2}-\frac{\epsilon}{12}}},
\end{eqnarray*}
outside a set of probability  at most
$
 C\left(-C^{'}\frac{R^{2}}{T^{\frac{3-\epsilon}{24}}\|\phi\|_{H^{s}}^{2}}\right).
$

\noindent When $|\xi_{1}|\sim |\xi_{6}|, \frac{|\xi_{7}|}{80}\leq |\xi_{1}|\leq80|\xi_{7}|$,
$|\xi_{7}|\leq 80|\xi_{8}|,$ by using  the  H\"older  inequality, (\ref{2.06}), (\ref{2.014}),
Lemmas 2.2,  we have that
\begin{eqnarray*}
&&I_{8}\leq CN_{6}^{-\frac{1}{20}}
\left(\prod_{j=1}^{6}\|J^{s}z_{j}\|_{L_{xt}^{\frac{48}{5-\epsilon}}}\right)\left(\prod_{j=7}^{8}\|J^{\sigma}
v_{j}\|_{L_{xt}^{8}}\right)
\|h\|_{L_{xt}^{\frac{8}{1+\epsilon}}}\nonumber\\&&\leq  CN_{6}^{-\frac{1}{20}}
\left(\prod_{j=1}^{6}\|J^{s}z_{j}\|_{L_{xt}^{\frac{48}{5-\epsilon}}}\right)
\left(\prod_{j=7}^{8}
\|v_{j}\|_{X_{\sigma,b}}\right)\|h\|_{X_{0,\frac{1}{2}-\frac{\epsilon}{12}}}\nonumber\\
&&\leq CN_{6}^{-\frac{1}{20}}R^{6}\left(\prod_{j=7}^{8}
\|v_{j}\|_{X_{\sigma,b}}\right)\|h\|_{X_{0,\frac{1}{2}-\frac{\epsilon}{12}}},
\end{eqnarray*}
outside a set of probability at most
$
 C{\rm exp}\left(-C^{'}\frac{R^{2}}{T^{\frac{5-\epsilon}{24}}\|\phi\|_{H^{s}}^{2}}\right).
$

\noindent When $\frac{|\xi_{7}|}{80}\geq|\xi_{1}|\sim |\xi_{6}|$,  this case can
be proved similarly to  case $\frac{|\xi_{7}|}{80}\geq|\xi_{1}|\sim |\xi_{2}|\geq 80|\xi_{3}|$
of Lemma 3.8.

\noindent When $|\xi_{1}|\sim |\xi_{7}|\geq80|\xi_{8}|$, this case can be proved
similarly to case  $|\xi_{1}|\sim |\xi_{6}|, \frac{|\xi_{7}|}{80}\leq |\xi_{1}|\leq80|\xi_{7}|$,
 $|\xi_{7}|\geq 80|\xi_{8}|$  of Lemma 3.8.

\noindent
When $|\xi_{1}|\sim|\xi_{8}|$, this case can be proved similarly to case
$|\xi_{1}|\sim |\xi_{6}|, \frac{|\xi_{7}|}{80}\leq |\xi_{6}|\leq80|\xi_{7}|$.

This completes the proof of Lemma 3.8.

\begin{Lemma}\label{Lemma3.9}
Let $s\geq\frac{17}{112}+\epsilon$ and $\sigma=\frac{3}{14}+2\epsilon$
and $z_{j}=\eta_{T}S(t)\phi^{\omega}(1\leq j\leq7,j\in N)$, $v_{8}=\eta v$. Then, we have that
\begin{eqnarray*}
\left|\int_{\SR}\int_{\SR}J^{\sigma}\partial_{x}
\left(\left(\prod_{j=1}^{7}z_{j}\right)v_{8}\right)h(x,t)dxdt\right|
\leq CT^{-\frac{3\epsilon}{100}}R^{7}\|v_{8}\|_{X_{\sigma,b}}\|h\|_{X_{0,\frac{1}{2}-\frac{\epsilon}{12}}},
\end{eqnarray*}
outside a set of probability at most
$
 C{\rm exp}\left(-C^{'}\frac{R^{2}}{\|\phi\|_{H^{s}}^{2}}\right).
$
\end{Lemma}
\noindent {\bf Proof.} We dyadically decompose $z_{j}$ with $(1\leq j\leq7,j\in N)$
 and $v_{8}$  and $h$ such that frequency supports are
$\left\{|\xi_{j}|\sim N_{j}\right\}$ for some dyadically $N_{j}\geq 1$
and we still denote them by $z_{j}$ with $(1\leq j\leq 7,j\in N)$, $v_{8}$ and $h$.
In this case, we divide the frequency into
$$|\xi_{l}|\geq {\rm max}\left\{|\xi_{j}|,1\leq j\leq 7,j\neq l,j\in N\right\}(1\leq l\leq 7,l\in N).$$
In this case, we divide the frequency  into
$|\xi_{l}|\geq {\rm max}\left\{|\xi_{j}|,1\leq j\leq 7,j\neq l,j\in N\right\}(1\leq l\leq 7,l\in N).$
Without loss of generality, we can assume that
$|\xi_{1}|\geq |\xi_{2}|\geq|\xi_{3}\geq |\xi_{5}|\geq|\xi_{6}|\geq|\xi_{7}|$ .

\noindent We define $$I_{9}=\left|\int_{\SR}\int_{\SR}J^{\sigma}\partial_{x}
\left(\left(\prod_{j=1}^{7}z_{j}\right)v_{8}\right)h(x,t)dxdt\right|.$$
\noindent When $|\xi_{1}|\geq80 |\xi_{2}|,$ $|\xi_{1}|\geq 80|\xi_{8}|,$
by using the  H\"older inequality, (\ref{2.04})-(\ref{2.05}), (\ref{2.016}),  Lemmas 2.4, 2.1, 2.3,
we obtain that
\begin{eqnarray*}
&&I_{9}\leq CN_{1}^{-\frac{1}{10}}
\left\|I^{1/2}I^{1/2}_{-}(J ^{s}z_{1},J ^{s}z_{2})\right\|_{L_{xt}^{2}}
\left(\prod_{j=4}^{7}\|J^{s}z_{j}\|_{L_{xt}^{\infty}}\right)
\|J^{-\frac{2}{7}}v_{8}\|_{L_{xt}^{\infty}}
\|I^{\frac{1-\epsilon}{2}}I^{\frac{1-\epsilon}{2}}_{-}(J ^{s}z_{3},h)\|_{L_{xt}^{2}}\nonumber\\&&\leq
CN_{1}^{-\frac{1}{10}}\left(\prod_{j=1}^{3}\|z_{j}\|_{X_{s,c}}\right)
\left(\prod_{j=4}^{7}\|J^{s}z_{j}\|_{L_{xt}^{\infty}}\right)\|v_{8}\|_{X_{\sigma,b}}
\|h\|_{X_{0,\frac{1}{2}-\frac{\epsilon}{12}}}\nonumber\\
&&\leq CN_{1}^{-\frac{1}{10}}T^{-\frac{3\epsilon}{100}}\left(\prod_{j=1}^{3}\|P_{N_{j}}\phi^{\omega}\|_{H^{s}}\right)
\left(\prod_{j=4}^{7}\|J^{s}z_{j}\|_{L_{xt}^{\infty}}\right)\|v_{8}\|_{X_{\sigma,b}}
\|h\|_{X_{0,\frac{1}{2}-\frac{\epsilon}{12}}}\nonumber\\
&&\leq CT^{-\frac{3\epsilon}{100}}R^{7}\|v_{8}\|_{X_{\sigma,b}}\|h\|_{X_{0,\frac{1}{2}-\frac{\epsilon}{12}}},
\end{eqnarray*}
outside a set of probability at most
$
C{\rm exp}\left(-C^{'}\frac{R^{2}}{\|\phi\|_{H^{s}}^{2}}\right).
$

\noindent When $|\xi_{1}|\geq80 |\xi_{2}|,$ $\frac{|\xi_{8}|}{80}\leq|\xi_{1}|\leq 80|\xi_{8}|,$
by using the H\"older inequality,  (\ref{2.05})-(\ref{2.06}),
(\ref{2.014}), Lemmas 2.4, 2.1, 2.2,  we have that
\begin{eqnarray*}
&&I_{9}\leq CN_{1}^{-\frac{1}{14}}
\left\|I^{1/2}I^{1/2}_{-}(J ^{s}z_{1},J ^{s}z_{2})\right\|_{L_{xt}^{2}}
\left(\prod_{j=3}^{7}\|z_{j}\|_{L_{xt}^{\frac{40}{2-\epsilon}}}\right)
\|J^{\sigma}v_{8}\|_{L_{xt}^{8}}
\|h\|_{L_{xt}^{\frac{8}{1+\epsilon}}}\nonumber\\&&\leq
CN_{1}^{-\frac{1}{14}}\left(\prod_{j=1}^{2}\|z_{j}\|_{X_{s,c}}\right)\left(\prod_{j=3}^{7}
\|z_{j}\|_{L_{xt}^{\frac{40}{2-\epsilon}}}\right)\|v_{8}\|_{X_{\sigma,b}}
\|h\|_{X_{0,\frac{1}{2}-\frac{\epsilon}{12}}}\nonumber\\
&&\leq CN_{1}^{-\frac{1}{14}}T^{-\frac{\epsilon}{50}}
\left(\prod_{j=1}^{2}\|P_{N_{j}}\phi^{\omega}\|_{H^{s}}\right)\left(\prod_{j=3}^{7}
\|z_{j}\|_{L_{xt}^{\frac{320}{17-7\epsilon}}}\right)\|v_{8}\|_{X_{\sigma,b}}
\|h\|_{X_{0,\frac{1}{2}-\frac{\epsilon}{12}}}\nonumber\\
&&\leq CT^{-\frac{\epsilon}{50}}R^{7}\|v_{8}\|_{X_{\sigma,b}}
\|h\|_{X_{0,\frac{1}{2}-\frac{\epsilon}{12}}},
\end{eqnarray*}
outside a set of probability at most
$
 C{\rm exp}\left(-C^{'}\frac{R^{2}}{\|\phi\|_{H^{s}}^{2}}\right).
$

\noindent When $\frac{|\xi_{8}|}{80}\geq|\xi_{1}|\geq80 |\xi_{2}|,$
by using the  H\"older  inequality,  (\ref{2.05})-(\ref{2.06}), Lemmas 2.4, 2.1, 2.2,
we have that
\begin{eqnarray*}
&&I_{9}\leq C
\left\|I^{1/2}I^{1/2}_{-}(J ^{\sigma}v_{8},z_{1})\right\|_{L_{xt}^{2}}
\left(\prod_{j=2}^{7}\|z_{j}\|_{L_{xt}^{\frac{48}{3-\epsilon}}}\right)
\|h\|_{L_{xt}^{\frac{8}{1+\epsilon}}}\nonumber\\&&\leq
C\|z_{1}\|_{X_{0,c}}
\left(\prod_{j=2}^{7}\|z_{j}\|_{L_{xt}^{\frac{48}{3-\epsilon}}}\right)
\|v_{8}\|_{X_{\sigma,b}}
\|h\|_{X_{0,\frac{1}{2}-\frac{\epsilon}{12}}}\nonumber\\
&&\leq CT^{-\frac{\epsilon}{100}}\|P_{N_{1}}\phi^{\omega}\|_{L^{2}}
\left(\prod_{j=2}^{7}\|z_{j}\|_{L_{xt}^{\frac{48}{3-\epsilon}}}\right)\|v_{8}\|_{X_{\sigma,b}}
\|h\|_{X_{0,\frac{1}{2}-\frac{\epsilon}{12}}}\nonumber\\
&&\leq CT^{-\frac{\epsilon}{100}}R^{7}\|v_{8}\|_{X_{\sigma,b}}
\|h\|_{X_{0,\frac{1}{2}-\frac{\epsilon}{12}}},
\end{eqnarray*}
outside a set of probability at most
$
 C{\rm exp}\left(-C^{'}\frac{R^{2}}{\|\phi\|_{H^{\epsilon}}^{2}}\right).
$

\noindent When $|\xi_{1}|\sim |\xi_{2}|\geq 80|\xi_{3}|,$ $|\xi_{1}|\geq \frac{|\xi_{8}|}{80},$
by using the H\"older  inequality, (\ref{2.05})-(\ref{2.07}),   Lemmas 2.4, 2.1, 2.2,
we get that
\begin{eqnarray*}
&&I_{9}\leq CN_{1}^{-\frac{1}{10}}
\|J^{s}z_{1}\|_{L_{xt}^{\frac{280}{17+7\epsilon}}}
\left\|I^{1/2}I^{1/2}_{-}(J^{s}z_{2},J^{s}z_{3})\right\|_{L_{xt}^{2}}
\left(\prod_{j=4}^{7}\|J^{s}z_{j}\|_{L_{xt}^{\frac{280}{17+7\epsilon}}}\right)
\|v_{8}\|_{L_{xt}^{\frac{28}{2-7\epsilon}}}
\|h\|_{L_{xt}^{\frac{8}{1+\epsilon}}}\nonumber\\&&\leq
CN_{1}^{-\frac{1}{10}}\|J^{s}z_{1}\|_{L_{xt}^{\frac{280}{17+7\epsilon}}}
\left(\prod_{j=2}^{3}\|z_{j}\|_{X_{s,c}}\right)
\left(\prod_{j=4}^{7}\|J^{s}z_{j}\|_{L_{xt}^{\frac{280}{17+7\epsilon}}}\right)\|v_{8}\|_{X_{\sigma,b}}
\|h\|_{X_{0,\frac{1}{2}-\frac{\epsilon}{12}}}\nonumber\\
&&\leq CN_{1}^{-\frac{1}{10}}T^{-\frac{\epsilon}{50}}\|J^{s}z_{1}\|_{L_{xt}^{\frac{280}{17+7\epsilon}}}
\left(\prod_{j=2}^{3}\|P_{N_{j}}\phi^{\omega}\|_{H^{s}}\right)
\left(\prod_{j=4}^{7}\|J^{s}z_{j}\|_{L_{xt}^{\frac{280}{17+7\epsilon}}}\right)\|v_{8}\|_{X_{\sigma,b}}
\|h\|_{X_{0,\frac{1}{2}-\frac{\epsilon}{12}}}\nonumber\\
&&\leq CT^{-\frac{\epsilon}{50}}R^{7}\|v_{8}\|_{X_{\sigma,b}}
\|h\|_{X_{0,\frac{1}{2}-\frac{\epsilon}{12}}},
\end{eqnarray*}
outside a set of probability at most
$
C{\rm exp}\left(-C^{'}\frac{R^{2}}{\|\phi\|_{H^{s}}^{2}}\right).
$

\noindent When $\frac{|\xi_{8}|}{80}\geq|\xi_{1}|\sim |\xi_{2}|\geq 80|\xi_{3}|,$
this case can be proved similarly to case $\frac{|\xi_{8}|}{80}\geq|\xi_{1}|\geq 80|\xi_{2}|$ of Lemma 3.9.

\noindent When $|\xi_{1}|\sim |\xi_{3}|\geq 80|\xi_{4}|,$ $|\xi_{1}|\geq \frac{|\xi_{8}|}{80},$
this case can be proved similarly to case $|\xi_{1}|\sim |\xi_{2}|\geq 80|\xi_{3}|,$
$|\xi_{1}|\geq \frac{|\xi_{8}|}{80}$ of Lemma 3.9.

\noindent When $\frac{|\xi_{8}|}{80}\geq|\xi_{1}|\sim |\xi_{3}|\geq 80|\xi_{4}|,$
this case can be proved similarly to  case $\frac{|\xi_{8}|}{80}\geq|\xi_{1}|\geq 80|\xi_{2}|$ of Lemma 3.9.

\noindent When $|\xi_{1}|\sim |\xi_{4}|\geq 80|\xi_{5}|,$ $|\xi_{1}|\geq \frac{|\xi_{8}|}{80}$,
 this case can be proved similarly to case $|\xi_{1}|\sim |\xi_{2}|\geq 80|\xi_{3}|,$
 $|\xi_{1}|\geq \frac{|\xi_{8}|}{80}$ of Lemma 3.9.

\noindent When $\frac{|\xi_{8}|}{80}\geq|\xi_{1}|\sim |\xi_{4}|\geq 80|\xi_{5}|,$
this case can be proved similarly to case $\frac{|\xi_{8}|}{80}\geq|\xi_{1}|\geq 80|\xi_{2}|$ of Lemma 3.9.

\noindent When $|\xi_{1}|\sim |\xi_{5}|\geq80|\xi_{6}|,|\xi_{5}|\geq\frac{|\xi_{8}|}{80}$,
 this case can be proved similarly to case $|\xi_{1}|\sim |\xi_{2}|\geq 80|\xi_{3}|,$
  $|\xi_{1}|\geq \frac{|\xi_{8}|}{80}$ of Lemma 3.9.

\noindent When $\frac{|\xi_{8}|}{80}\geq|\xi_{1}|\sim |\xi_{5}|\geq80|\xi_{6}|$,
this case can be proved similarly to  case $\frac{|\xi_{8}|}{80}\geq|\xi_{1}|\geq 80|\xi_{2}|$ of Lemma 3.9.

\noindent When $|\xi_{1}|\sim |\xi_{6}|\geq80|\xi_{7}|, |\xi_{1}|\geq\frac{|\xi_{8}|}{80}$,
this case can be proved similarly to case $|\xi_{1}|\sim |\xi_{2}|\geq 80|\xi_{3}|,$
 $|\xi_{1}|\geq \frac{|\xi_{8}|}{80}$ of Lemma 3.9.

\noindent When $\frac{|\xi_{8}|}{80}\geq|\xi_{1}|\sim |\xi_{6}|\geq80|\xi_{7}|$,
this case can be proved similarly to  case  $\frac{|\xi_{8}|}{80}\geq|\xi_{1}|\geq 80|\xi_{2}|$ of Lemma 3.9.

\noindent When $|\xi_{1}|\sim |\xi_{7}|, |\xi_{1}|\geq80|\xi_{8}|$,  by using  the  H\"older  inequality,
 (\ref{2.05})-(\ref{2.06}),   Lemmas 2.4, 2.1, 2.2, we get that
\begin{eqnarray*}
&&I_{9}
\leq CN_{7}^{-\frac{1}{2}}\left(\prod_{j=2}^{7}\|J^{s}z_{j}\|_{L_{xt}^{\frac{48}{3-\epsilon}}}\right)
\left\|I^{1/2}I^{1/2}_{-}(J^{s}z_{1},J^{\sigma}v_{8})\right\|_{L_{xt}^{2}}
\|h\|_{L_{xt}^{\frac{8}{1+\epsilon}}}\nonumber\\&&\leq CN_{7}^{-\frac{1}{2}}
\|z_{1}\|_{X_{s,c}}\left(\prod_{j=2}^{7}\|J^{s}z_{j}\|_{L_{xt}^{\frac{48}{3-\epsilon}}}\right)
\|v_{8}\|_{X_{\sigma,b}}\|h\|_{X_{0,\frac{1}{2}-\frac{\epsilon}{12}}}\nonumber\\
&&\leq CN_{7}^{-\frac{1}{2}}T^{-\frac{\epsilon}{100}}\|P_{N_{1}}\phi^{\omega}\|_{H^{s}}
\left(\prod_{j=2}^{7}\|J^{s}z_{j}\|_{L_{xt}^{\frac{48}{3-\epsilon}}}\right)
\|v_{8}\|_{X_{\sigma,b}}\|h\|_{X_{0,\frac{1}{2}-\frac{\epsilon}{12}}}\nonumber\\
&&\leq CT^{-\frac{\epsilon}{100}}R^{7}\|v_{8}\|_{X_{\sigma,b}}\|h\|_{X_{0,\frac{1}{2}-\frac{\epsilon}{12}}},
\end{eqnarray*}
outside a set of probability at most
$ C\left(-C^{'}\frac{R^{2}}{\|\phi\|_{H^{s}}^{2}}\right).
$

\noindent
When $|\xi_{1}|\sim|\xi_{7}|,\frac{|\xi_{8}|}{80}\leq |\xi_{1}|\leq80|\xi_{8}|$, by using  the  H\"older   inequality, (\ref{2.06})
and (\ref{2.015}), Lemma 2.2,   we have that
\begin{eqnarray*}
&&I_{9}\leq CN_{1}^{-\frac{1}{20}}
\left(\prod_{j=1}^{7}\|J^{s}z_{j}\|_{L_{xt}^{\frac{56}{6-\epsilon}}}\right)\|J^{\sigma}
v_{8}\|_{L_{xt}^{8}}
\|h\|_{L_{xt}^{\frac{8}{1+\epsilon}}}\nonumber\\&&\leq  CN_{1}^{-\frac{1}{20}}
\left(\prod_{j=1}^{7}\|J^{s}z_{j}\|_{L_{xt}^{\frac{56}{6-\epsilon}}}\right)
\|v_{8}\|_{X_{\sigma,b}}\|h\|_{X_{0,\frac{1}{2}-\frac{\epsilon}{12}}}\leq
CR^{7}\|v_{8}\|_{X_{\sigma,b}}\|h\|_{X_{0,\frac{1}{2}-\frac{\epsilon}{12}}},
\end{eqnarray*}
outside a set of probability at most
$
 C\left(-C^{'}\frac{R^{2}}{T^{\frac{6-\epsilon}{28}}\|\phi\|_{H^{s}}^{2}}\right).
$

\noindent
When $\frac{|\xi_{8}|}{80}\geq|\xi_{1}|\sim|\xi_{7}|$,
this case can be proved similarly to   $\frac{|\xi_{8}|}{80}\geq|\xi_{1}|\geq 80|\xi_{2}|$ of Lemma 3.9.

\noindent When $|\xi_{1}|\sim|\xi_{8}|$, this case can be proved similarly to case $|\xi_{1}|\sim|\xi_{7}|,
\frac{|\xi_{8}|}{80}\leq |\xi_{1}|\leq80|\xi_{8}|$ of Lemma 3.9.

This completes the proof of Lemma 3.9.
\bigskip

\begin{Lemma}\label{Lemma3.10}
Let  $\sigma=\frac{3}{14}+2\epsilon$ and $0<T<1$. Then, we have that
\begin{eqnarray}
\left\|(\eta v+\eta_{T}z)
\partial_{x}\left[(\eta v+\eta_{T} z)^{7}\right]\right\|_{X_{\sigma,-\frac{1}{2}+\frac{\epsilon}{12}}}
\leq CT^{-\frac{3\epsilon}{100}}\left[\sum_{m=0}^{8}\|v\|_{X_{\sigma,b}}^{m}R^{8-m}\right]\label{3.01}
\end{eqnarray}
outside a set of probability at most
$
C{\rm exp}\left(-C^{'}\frac{R^{2}}{\|\phi\|_{H^{s}}^{2}}\right).
$
\end{Lemma}
\noindent{\bf Proof.} By duality, to prove (\ref{3.01}), it suffices to prove that
\begin{eqnarray*}
\left|\int_{\SR}\int_{\SR}(\eta v+\eta_{T}z)
\partial_{x}\left[(\eta v+\eta_{T} z)^{7}\right]hdxdt\right|
\leq CT^{-\frac{3\epsilon}{100}}\left[\sum_{m=0}^{8}\|v\|_{X_{\sigma,b}}^{m}R^{8-m}\right]
\|h\|_{X_{0,\frac{1}{2}-\frac{\epsilon}{12}}},
\end{eqnarray*}
which can be obtained from Lemmas 3.1-3.9.

We have completed the proof of Lemma 3.10.

\noindent {\large\bf 4. Proof of Theorem  1.1}

\setcounter{equation}{0}

 \setcounter{Theorem}{0}

\setcounter{Lemma}{0}

\setcounter{section}{4}
In this section, we are in a position to prove Theorem 1.1.

\noindent {\bf Proof of  Theorem 1.1.} Assume that $\phi^{\omega}$ is the randomization of $\phi$,
 which satisfies (\ref{1.02}) and belongs to $H^{s}(\R)$ almost surely.
Now we consider the Cauchy problem for (\ref{1.01}) with $u(x,0)=\phi^{\omega}$.
Let $z(t)=z^{\omega}(t)=S(t)\phi^{\omega}$ and $v(t)=u(t)-z(t)$, (\ref{1.01}) can be rewritten as follows:
\begin{eqnarray}
&&v_{t}+v_{xxx}+(v+z)\partial_{x}\left[(v+z)^{7}\right]=0,\label{4.01}\\
&&v(x,0)=0\label{4.02}.
\end{eqnarray}
Note that (\ref{4.01})-(\ref{4.02}) are equivalent to the following integral equation
\begin{eqnarray}
v(t)=\int_{0}^{t}S(t-\tau)(v+z)\partial_{x}\left[(v+z)^{7}\right]d\tau.\label{4.03}
\end{eqnarray}
Obviously, $v$ satisfies (\ref{4.03}) on $[-T,T]$ if v satisfies
\begin{eqnarray}
v(t)=\eta_{T}(t)\int_{0}^{t}S(t-\tau)(\eta v+\eta_{T}z)
\partial_{x}\left[(\eta v+\eta_{T} z)^{7}\right]d\tau.\label{4.04}
\end{eqnarray}
for some $T\ll1$.
We define
\begin{eqnarray}
\Gamma v=\eta_{T}(t)\int_{0}^{t}S(t-\tau)(\eta v+\eta_{T}z)
\partial_{x}\left[(\eta v+\eta_{T} z)^{7}\right]d\tau.\label{4.05}
\end{eqnarray}
By using Lemmas 2.4, 3.10 and the Young inequality, we get that
\begin{eqnarray}
&&\|\Gamma v\|_{X^{\sigma,b}}\leq CT^{\frac{\epsilon}{24}}
\left\|(\eta v+\eta_{T}z)\partial_{x}\left[(\eta v+\eta_{T}z)^{7}\right]
\right\|_{X_{\sigma,-\frac{1}{2}+\frac{\epsilon}{12}}}\nonumber\\&&\leq CT^{\frac{\epsilon}{24}}
\left[\sum_{m=0}^{8}\|v\|_{X_{\sigma,b}}^{m}R^{8-m}\right]\nonumber\\&&
\leq CT^{\frac{\epsilon}{24}}T^{-\frac{3\epsilon}{100}}\left[\sum_{m=1}^{7}\|v\|_{X_{\sigma,b}}^{m}R^{8-m}
+\|v\|_{X_{\sigma,b}}^{8}+R^{8}\right]\nonumber\\
&&\leq CT^{\frac{7\epsilon}{600}}\left[7\left(\frac{m}{8}\|v\|_{X_{\sigma,b}}^{8}
+\frac{8-m}{8}R^{8}\right)+\|v\|_{X_{\sigma,b}}^{8}+R^{8}\right]\nonumber\\&&
\leq CT^{\frac{7\epsilon}{600}}
\left(\|v\|_{X_{\sigma,b}}^{8}+R^{8}\right),\label{4.06}
\end{eqnarray}
outside a set of probability at most
$
C{\rm exp}\left(-C^{'}\frac{R^{2}}{\|\phi\|_{H^{s}}^{2}}\right).
$
Similarly,  by using 2.4, 3.10 and the Young inequality,  we conclude that
\begin{eqnarray}
\|\Gamma v_{1}-\Gamma v_{2}\|_{X_{\sigma,b}}\leq C
T^{\frac{7\epsilon}{600}}\|v_{1}-v_{2}\|_{X_{\sigma,b}}
\left(\|v_{1}\|_{X_{\sigma,b}}^{7}+\|v_{2}\|_{X_{\sigma,b}}^{7}+R^{7}\right),\label{4.07}
\end{eqnarray}
outside a set of probability at most $C{\rm exp}\left(-c\frac{R^{2}}{\|\phi\|_{H^{s}}^{2}}\right)$.
Let $B_{1}=\left\{u\mid\|u\|_{X_{\sigma,b}}\leq1\right\}$, for $T\ll1$, we choose $R=R(T)\sim T^{-\frac{7\epsilon}{4800}}$
such that
\begin{eqnarray}
&&CT^{\frac{7\epsilon}{600}}(1+R^{8})\leq 1,CT^{\frac{7\epsilon}{600}}(2+R^{7})
\leq \frac{1}{2}\label{4.08}.
\end{eqnarray}
Thus, for $v,v_{1},v_{2}\in B_{1}$, combining (\ref{4.06}), (\ref{4.07})  with  (\ref{4.08}), we have that
\begin{eqnarray}
&&\|\Gamma v\|_{X_{\sigma,b}}\leq1,\|\Gamma v_{1}-\Gamma v_{2}\|_{X_{\sigma,b}}
\leq\frac{1}{2}\|v_{1}-v_{2}\|_{X_{\sigma,b}}\label{4.09}.
\end{eqnarray}
outside an exceptional set of probability at most
\begin{eqnarray}
C{\rm exp}\left(-C^{'}\frac{R^{2}}{\|\phi\|_{H^{s}}^{2}}\right)\sim C{\rm exp}
\left(-\frac{C^{'}}{T^{\frac{7\epsilon}{2400}}\|\phi\|_{H^{s}}^{2}}\right)
\end{eqnarray}
Consequently, let $\Omega_{T}$  be the complement of the exceptional set,
for $\omega \in \Omega_{T}$, there exists a unique $v^{\omega}\in B_{1}$ such
that $\Gamma v^{\omega}=v^{\omega}$.

This completes the proof of Theorem 1.1.

\bigskip
\leftline{\large \bf Acknowledgments}

\noindent This work is supported by NSF  grant  1620449 and the Natural Science Foundation of China
  grant 11401180 and 11371167. The first author is also
 supported by   the Young Core Teachers program of Henan Normal University and  15A110033.

\end{document}